\documentclass{amsart}
\usepackage{amsmath} 
\usepackage{amssymb}
\usepackage{mathrsfs}

\newtheorem{theorem}{Theorem}[section] 
\newtheorem{claim}[theorem]{Claim}

\newtheorem{conclusion}[theorem]{Conclusion}
\newtheorem{observation}[theorem]{Observation}

\theoremstyle{definition}
\newtheorem{definition}[theorem]{Definition}

\theoremstyle{remark}
\newtheorem{remark}[theorem]{Remark}
\newtheorem{question}[theorem]{Question}
\newtheorem{notation}[theorem]{Notation}

\newcommand{\up}{{\rm up}}
\newcommand{\uf}{{\rm uf}}
\newcommand{\md}{{\rm md}}
\newcommand{\dn}{{\rm dn}}
\newcommand{\lt}{{\rm lt}}

\newcommand{\acc}{{\rm acc}}

\newcommand{\Qr}{{\rm Qr}}

\newcommand{\Col}{{\rm Col}}
\newcommand{\otp}{{\rm otp}}
\newcommand{\Max}{{\rm Max}}

\newcommand{\Ord}{{\rm Ord}}

\newcommand{\range}{{\rm range}}

\newcommand{\nacc}{{\rm nacc}}

\newcommand{\Per}{{\rm Pr}}

\newcommand{\Dom}{{\rm Dom}}
\newcommand{\Rang}{{\rm Rang}}

\newcommand{\bfi}{{\bold i}}

\newcommand{\rest}{{\restriction}}
\newcommand{\dom}{{\rm dom}}

\newcommand{\wilog}{{\rm without loss of generality}}
\newcommand{\Wilog}{{\rm Without loss of generality}}

\newcommand{\then}{{\underline{then}}}
\newcommand{\when}{{\underline{when}}}
\newcommand{\Then}{{\underline{Then}}}

\newcommand{\Iff}{{\underline{iff}}}
\newcommand{\mn}{{\medskip\noindent}}
\newcommand{\sn}{{\smallskip\noindent}}

\newcommand{\cF}{{\mathscr F}}
\newcommand{\cG}{{\mathscr G}}

\newcommand{\bbN}{{\mathbb N}}

\newcommand{\cP}{{\mathscr P}}

\newcommand{\varp}{{\varepsilon}} 
\newcommand{\cU}{{\mathscr U}}

\newcommand{\cW}{{\mathscr W}}

\newcommand{\cf}{{\rm cf}}
\newcommand{\pr}{{\rm pr}}

\newcount\skewfactor
\def\mathunderaccent#1#2 {\let\theaccent#1\skewfactor#2
\mathpalette\putaccentunder}
\def\putaccentunder#1#2{\oalign{$#1#2$\crcr\hidewidth
\vbox to.2ex{\hbox{$#1\skew\skewfactor\theaccent{}$}\vss}\hidewidth}}

\newenvironment{PROOF}[2][\proofname.]
   {\begin{proof}[#1]\renewcommand{\qedsymbol}{\textsquare$_{\rm #2}$}}
   {\end{proof}}

\begin{document}

\title {The colouring existence theorem revisited}
\author {Saharon Shelah}
\address{Einstein Institute of Mathematics\\
Edmond J. Safra Campus, Givat Ram\\
The Hebrew University of Jerusalem\\
Jerusalem, 91904, Israel\\
 and \\
 Department of Mathematics\\
 Hill Center - Busch Campus \\ 
 Rutgers, The State University of New Jersey \\
 110 Frelinghuysen Road \\
 Piscataway, NJ 08854-8019 USA}
\email{shelah@math.huji.ac.il}
\urladdr{http://shelah.logic.at}
\thanks{The author would like to thank the Israel Science Foundation
  for partial support of this research.
The author thanks Alice Leonhardt for the beautiful typing.
  First typed June 13, 2013. Paper 1027}





\date{June 29, 2018}

\begin{abstract}
We prove a better colouring theorem for $\aleph_4$ and even
$\aleph_3$.  This has a general topology consequence.
\end{abstract}

\maketitle
\numberwithin{equation}{section}
\setcounter{section}{-1}
\newpage

\section {Introduction}
\bigskip

\subsection {Background}\
\bigskip

Our aim is to improve some colouring theorems of \cite{Sh:327},
\cite[Ch.III,\S4]{Sh:g}, they continue Tordor\'cevi\'c \cite{To2}
(introducing the walks) and \cite{Sh:280}, \cite[\S3]{Sh:276} (and 
\cite{Sh:572}), see history in
\cite{Sh:g}, \cite[\S10]{Sh:E12}.  After these works Moore \cite{Mo06} prove
$\aleph_1 \mapsto [\aleph_1;\aleph_1]^2_{\aleph_0}$; Eisworth \cite{Eis13} and
Rinot \cite{Rin13} prove equivalence of some colouring theorems on successor of
singular cardinals.

Our aim is to prove better colouring theorems on successor of regular cardinals
(when not too small),
e.g. $\Pr_1(\aleph_3,\aleph_3,\aleph_3,(\aleph_0,\aleph_1))$, see
\S1.  We have looked at the matter again because Juhasz-Shelah
\cite{JuSh:1025} need such theorem in order to solve a problem in
general topology, see \ref{d47}(3).
 
\bigskip

\subsection {Results}\
\bigskip

The paper is self contained.

Here we formulate $\Pr_\ell(\lambda,\mu,\sigma,\bar\theta)$ where
$\bar\theta$ is a pair $(\theta_0,\theta_1)$ of cardinals rather than a single
cardinal $\theta$ and prove
e.g. $\Pr(\lambda,\lambda,\lambda,(\theta,\theta^+))$ when $\lambda =
\theta^{+3}$ and $\theta$ is regular.

That is, we shall prove (see Definition \ref{d32} and Conclusion \ref{d47}(1)):

\begin{theorem}  
\label{y2}
1) For any regular $\kappa$ we have
$\Pr_1(\kappa^{+4},\kappa^{+4},\kappa^{+4},\kappa^+)$.

\noindent
2) For any regular $\kappa$ we have
($\Pr_1(\kappa^{+4},\kappa^{+4},\kappa^{+4},(\kappa,\kappa^+))$ and) 
$\Pr_{0,0}(\kappa^{+4},\kappa^{+4},\kappa^{+4},(\aleph_0,\kappa^+))$. 
\end{theorem}

\begin{remark}
Note that the statement $\Pr_0(\kappa^{+4},\kappa^{+4},2,\kappa^+)$ is also 
called by Juhasz $\Col(\kappa^{+4},\kappa)$, see more in the end of \S1.
\end{remark}

\noindent
Moreover by \ref{d50} in \ref{y2}(2) we 
can replace $\kappa^{+4}$ by $\kappa^{+3}$,
(thus half solving Problem 1 of \cite{JuSh:1025}, i.e. for $\aleph_3$
though not for $\aleph_2$) so we naturally ask:
\begin{question}
1) Do we have $\Pr_1(\aleph_2,\aleph_2,\sigma,\aleph_1)$ for $\sigma =
   \aleph_2$?  For $\sigma = 2$?

\noindent
2) Do we have at least
$\Pr^{\uf}_{0,0}(\aleph_2,\aleph_2,2,(\aleph_0,\aleph_1))$?

Concerning the result of Juhasz-Shelah \cite{JuSh:1025} by using \ref{d43}(1)
instead of \cite[Ch.III,\S4]{Sh:g} we can deduce
$\Pr_0(\aleph_4,\aleph_4,2,(\aleph_0,\aleph_1))$ which is sufficient
for the topological result there.
Moreover by \ref{e16} + \ref{d35} even
$\Pr_0(\aleph_3,\aleph_3,2,(\aleph_0,\aleph_1))$ holds, see \ref{d47}
so there is a topological space as desired in \cite{JuSh:1025} with
weight $\aleph_3$, see \ref{d50}(2).

We can also generalize the other conclusion of \cite[Ch.III,\S4]{Sh:g}
replacing $\theta$ by $(\theta_0,\theta_1)$.  This may be dealt with
later.  Also in \cite{Sh:F1422} and better \cite{Sh:F1783} we intend 
to improve \ref{d50} for most cardinals.

We thank Shimoni Garti and the referee for pointing out many missing points. 
\end{question}
\newpage

\section {Definitions and some connections}

\begin{definition}  
\label{d32}
Assume $\lambda \ge \mu \ge \sigma + \theta_0 + \theta_1,\bar\theta
   = (\theta_0,\theta_1)$; if $\theta_0 = \theta_1$ we may write
   $\theta_0$ instead of $\bar\theta$.

\noindent
1) Let $\Pr_0(\lambda,\mu,\sigma,\bar\theta)$ mean that there is $\bold
c:[\lambda]^2 \rightarrow \sigma$ witnessing it which means:
\mn
\begin{enumerate}
\item[$(*)_{\bold c}$]  if (a) then (b) where:
\sn
\begin{enumerate}
\item[$(a)$]
\begin{enumerate}
\item[$(\alpha)$]   for $\iota = 0,1,\bar\zeta^\iota = \langle
  \zeta^\iota_{\alpha,i}:\alpha < \mu,i < \bold i_\iota\rangle$ is a
  sequence without repetitions of ordinals $< \lambda$ and
  $\Rang(\bar\zeta^0),\Rang(\bar\zeta^1)$ are disjoint 
 and $\bold i_0 < \theta_0,\bold i_1 < \theta_1$
\sn
\item[$(\beta)$]  $h:\bold i_0 \times \bold i_1 \rightarrow \sigma$
\end{enumerate}
\sn
\item[$(b)$]   for some $\alpha_0 < \alpha_1 < \mu$ we have:
\begin{itemize}
\item  if $i_0 < \bold i_0$ and $i_1 < \bold i_1$ then
$\bold c \{\zeta^0_{\alpha_0,i_0},\zeta^1_{\alpha_1,i_1}\} =
h(i_0,i_1)$.
\end{itemize}
\end{enumerate}
\end{enumerate}
\mn
2) For $\iota \in \{0,1\}$ let
$\Pr_{0,\iota}(\lambda,\mu,\sigma,\bar\theta)$ be defined similarly
but we replace $(a)(\beta)$ and $(b)$ by $(a)(\beta)'$ and $(b)'$, where
\mn
\begin{enumerate}
\item[$(a)$] 
\begin{enumerate}
\item[$(\beta)'$]  $h:\bold i_\iota \rightarrow \sigma$ 
\end{enumerate}
\sn
\item[$(b)'$]   for some $\alpha_0 < \alpha_1 < \mu$ we have
\begin{enumerate}
\item[$\bullet'$]   if $i_0 < \bold i_0$ and $i_1 < \bold i_1$ then
$\bold c \{\zeta^0_{\alpha_0,i_0},\zeta^1_{\alpha_1,i_1}\} = h(i_\iota)$.
\end{enumerate}
\end{enumerate}
\mn
3) Let $\Pr^{\uf}_{0,\iota}(\lambda,\mu,\sigma,\bar\theta)$ mean that some
$\bold c:[\lambda]^2 \rightarrow \sigma$ witnesses it which means:
\mn
\begin{enumerate}
\item[$(*)^{\uf}_{\bold c}$]  if (a) then (b) where
\mn
\begin{enumerate}
\item[$(a)$]   $(\alpha) \quad$ as above
\sn
\item[${{}}$]  $(\beta) \quad h:\bold i_\iota \rightarrow \sigma$ and
  $D$ is an ultrafilter on $\bold i_{1-\iota}$
\sn
\item[$(b)$]   for some $\alpha_0 < \alpha_1 < \mu$ we have
\sn
\item[${{}}$]  $\bullet \quad$ if $i < \bold i_\iota$ then $\{j <
  \bold i_{1-\iota}:\bold
  c\{\zeta^\iota_{\alpha_\iota,i},\zeta^{1-\iota}_{\alpha_{1-\iota,j}}\}
  = h(i)\}$ belongs to $D$.
\end{enumerate}
\end{enumerate}
\end{definition}

\begin{definition}
\label{d321}
Assume $\lambda \ge \mu \ge \sigma + \theta_0 + \theta_1,\bar\theta =
(\theta_0,\theta_1)$.  Let $\Pr_1(\lambda,\mu,\sigma,\bar\theta)$ mean
that there is $\bold c:[\lambda]^2 \rightarrow \sigma$ witnessing it,
which means:
\mn
\begin{enumerate}
\item[$(*)_{\bold c}$]  if (a) then (b), where:
\sn
\begin{enumerate}
\item[(a)]  for $\iota=0,1,\bold i_\iota < \theta_\iota$ and
  $\bar\zeta^\iota = \langle \zeta^\iota_{\alpha,i}:\alpha < \mu,i <
  \bold i_\iota\rangle$ are sequences of ordinals of $\lambda$ without
  repetitions, $\Rang(\bar\zeta^\iota)$ are disjoint and $\gamma < \sigma$
\sn
\item[(b)]  there are $\alpha_0 < \alpha_1 < \mu$ such that $\forall
  i_0 < \bold i_0,\forall i_1 < \bold i_1,\bold c
  \{\zeta^0_{\alpha_0,i_0},\zeta^1_{\alpha_1,i_1}\} = \gamma$.
\end{enumerate}
\end{enumerate}
\end{definition}

\begin{remark}  
\label{d34}
1) So if $\theta_0 = \theta = \theta_1$ and $\bar\theta =
(\theta_0,\theta_1)$ then for $\ell \in \{0,1\},\Pr_\ell(\lambda,\mu,
\sigma,\bar\theta)$ is $\Pr_\ell(\lambda,\mu,\sigma,\theta)$ from 
\cite[Ch.III]{Sh:g}.

\noindent
2) We do not write down the monotonicity and trivial implications
 concerning Definitions \ref{d32} and \ref{d35} below.

\noindent
3) The disjointness of $\{\zeta^0_{\alpha,i}:\alpha < \mu,i < \bold
i_0\},\{\zeta^1_{\alpha,i}:\alpha < \mu,i < \bold i_1\}$ in Definition
\ref{d32}(1),(a),$(\alpha)$ and \ref{d32}(2), \ref{d32}(3) and
\ref{d321}(a) is not really necessary. 
\end{remark}

\begin{notation}
$\pr:\Ord \times \Ord \rightarrow \Ord$ is the standard pairing
  function.
\end{notation}

\noindent
Variants are
\begin{definition}
\label{d35}
Let $\lambda \ge \mu \ge \sigma + \theta_0 + \theta_1$ and $\bar\theta
= (\theta_0,\theta_1)$.

\noindent
1) Let $\Qr_0(\lambda,\mu,\sigma,\bar\theta)$ mean that there is
   $\bold c:[\lambda]^2 \rightarrow \sigma$ witnessing it which means:
\mn
\begin{enumerate}
\item[$(*)_{\bold c}$]  if (a) then (b) where
\sn
\begin{enumerate}
\item[$(a)$]
\begin{enumerate}
\item[$(\alpha)$]  $u^\iota_\alpha \in [\lambda]^{<
  \theta_\iota}$ for $\iota < 2$ and $\alpha < \mu$
\sn
\item[$(\beta)$]  $u_\alpha = u^0_\alpha \cup u^1_\alpha$ for every
  $\alpha < \mu$
\sn
\item[$(\gamma)$]  $\langle u_\alpha:\alpha < \mu\rangle$
 are pairwise disjoint
\sn
\item[$(\delta)$]  $h^\iota_\alpha:u^\iota_\alpha
  \rightarrow \sigma$ for $\iota < 2,\alpha < \mu$ 
and $\pr:\sigma \times \sigma \rightarrow \sigma$
\end{enumerate}
\sn
\item[$(b)$]  for some $\alpha_0 < \alpha_1 < \mu$ for every
  $(\zeta_0,\zeta_1) \in (u^0_{\alpha_0} \times u^1_{\alpha_1})$ we have
 $\zeta_0 < \zeta_1$ and $\bold c\{\zeta_0,\zeta_1\} =
  \pr(h^0_{\alpha_0}(\zeta_0),h^1_{\alpha_1}(\zeta_1))$.
\end{enumerate}
\end{enumerate}
\mn
2) Let $\Qr_{0,\iota}(\lambda,\mu,\sigma,\bar\theta)$ be defined
similarly but each $h^{1-\iota}_\alpha$ is constant.

\noindent
3) Let $\Qr_1(\lambda,\mu,\sigma,\bar\theta)$ be defined as above but each
$h^0_\alpha$ and each $h^1_\alpha$ is a constant function.

\noindent
4) Let $\Qr^{\uf}_{0,\iota}(\lambda,\mu,\sigma,\bar\theta)$ be defined
parallely to Definition \ref{d32}. 
\end{definition}

\noindent
So, e.g.
\begin{observation}
\label{d36}
1) If $\cf(\mu) \ge \sigma^+$, \then \,
$\Pr_1(\lambda,\mu,\sigma,\bar\theta)$ is equivalent to
$\Qr_1(\lambda,\mu,\sigma,\bar\theta)$.

\noindent
2) Recall that $\Pr_\ell(\lambda,\mu,\sigma,\theta)$ is
 $\Pr_\ell(\lambda,\mu,\sigma,(\theta,\theta))$.

\noindent
3) $\Qr_0(\lambda,\mu,\sigma,\bar\theta)$ implies 
$\Pr_0(\lambda,\mu,\sigma,\bar\theta)$; similarly for the other
variants, $\Qr_{0,\iota},\Qr^{\uf}_{0,\iota}$.
\end{observation}

\begin{PROOF}{\ref{d36}}
Should be clear.
\end{PROOF}

\begin{observation}  
\label{d40}
Let $\bar\theta = (\theta_0,\theta_1)$ and $\iota \in \{0,1\}$.

\noindent
1) If $\iota < 2,\partial < \theta_\iota 
\Rightarrow \sigma^\partial < \cf(\mu)$ and $\theta_0,\theta_1 < 
\cf(\mu)$, \then \, 
$\Pr_{0,\iota}(\lambda,\mu,\sigma,\bar\theta)$ is equivalent to 
$\Qr_{0,\iota}(\lambda,\mu,\sigma,\bar\theta)$.

\noindent 
2) If $\partial < \theta_0 + \theta_1 \Rightarrow \sigma^\partial <
\cf(\mu)$, \then \, $\Pr_0(\lambda,\mu,\sigma,\bar\theta) \Leftrightarrow
\Qr_0(\lambda,\mu,\sigma,\bar\theta)$. 
\end{observation}

\begin{PROOF}{\ref{d40}}
Obvious but we elaborate.

\noindent
1) By \ref{d36}(3) we have one implication; so assume
$\Pr_{0,\iota}(\lambda,\mu,\sigma,\bar\theta)$ and we shall prove
$\Qr_{0,\iota}(\lambda,\mu,\sigma,\bar\theta)$, so let $u_\alpha =
u^0_\alpha \cup u^1_\alpha$ for $\alpha < \mu$ and
$h^\iota_d:u^\iota_d \rightarrow \sigma$ and $\pr:\sigma \times
\sigma$ be as in Definition \ref{d35}(1) and each $h^{1-\iota}_\alpha$
is constant.

We should prove that there are $\alpha_0 < \alpha_1 < \mu$ as promised
in Definition \ref{d35}(2).  As
$|u^{1-\varp}_\alpha| < \theta_{1-\iota}$ and $\theta_{1-\iota} <
\cf(\mu)$, \wilog \, for some $\varp_{1-\iota} < \theta_{1-\iota}$ we
have $\alpha < \mu \Rightarrow \otp(u^{1-\iota}_\alpha) =
\varp_{1-\iota}$.  As $\theta_\iota < \cf(\mu)$ hence
\wilog \, for some $\varp_\iota < \theta_\iota$ we have $\alpha < \mu
\Rightarrow \otp(u^\iota_\alpha) = \varp_\iota$.  Moreover, noting
$\sigma^{|\varp_\iota|} < \cf(\mu)$, \wilog \, $\{(\otp(\zeta \cap
u^\iota_\alpha),h_\alpha(\zeta)):\zeta \in u^\iota_\alpha\}$ is the
same for all $\alpha < \mu$.  Now we can apply
$\Pr_{0,\iota}(\lambda,\mu,\sigma,\bar\theta)$.

\noindent
2) Similarly.
\end{PROOF}

\begin{claim}  
\label{d43}
1) Let $\iota < 2$.  
If $\Pr_1(\lambda,\mu,\sigma_1,\bar\theta)$ and
 $\lambda = \mu = \cf(\mu),\bar\theta = (\theta_0,\theta_1),
\theta = \theta_0 + \theta_1 < \mu$ and 
$2^\chi \ge \lambda,\chi^{< \theta_\iota} + (\sigma_2)^{<
  \theta_\iota} \le \sigma_1$ and $\chi^{< \theta_\iota} < \mu$ and
$(\sigma_2)^{< \theta_\iota} < \mu$ \then \, 
$\Pr_{0,\iota}(\lambda,\mu,\sigma_2,\bar\theta)$ and
$\Qr_{0,\iota}(\lambda,\mu,\sigma_2,\theta)$.

\noindent
1A) If the assumptions of part (1) holds for both $\iota=0$ and
$\iota=1$, \then \, we can conclude
$\Pr_0(\lambda,\mu,\sigma_2,\bar\theta)$ and
$\Qr_0(\lambda,\mu,\sigma_2,\bar\theta)$. 

\noindent
2) If $\lambda = \sigma^+$ and $\sigma = \sigma^{< \theta_\iota}$ \then \,
$\Pr_{0,\iota}(\lambda,\lambda,\sigma,\bar\theta)$ implies
$\Pr_{0,\iota}(\lambda,\lambda,\lambda,\bar\theta)$.

\noindent
3) If $\lambda = \sigma^+$ and $\sigma = \sigma^{<(\theta_0 +
   \theta_1)}$ \then \, $\Pr_0(\lambda,\lambda,\sigma,\bar\theta)$ implies
$\Pr_0(\lambda,\lambda,\lambda,\bar\theta)$.

\noindent
4) If $\Pr_1(\lambda,\mu,\sigma,\bar\theta)$ and $\sigma \le \chi =
\chi^{<(\theta_0 + \theta_1)} < \lambda \le 2^\chi$ \then \,
$\Pr_0(\lambda,\mu,\sigma,\bar\theta)$.

\noindent
5) If $\Pr_1(\lambda,\lambda,\lambda,\bar\theta),\lambda = \partial^+$
and $\partial = \partial^{<(\theta_0 + \theta_1)}$ \then \, 
$\Pr_0(\lambda,\lambda,\lambda,\bar\theta)$. 
\end{claim}

\begin{remark}
\label{d44}
1) Claim \ref{d43}(1) is similar to \
\cite[Ch.III,4.5(3),pg.169-170]{Sh:g} but we shall elaborate.

\noindent
2) The condition $\lambda = \mu$ can be omitted if we systematically
use $\bold c:\lambda \times \lambda \rightarrow \sigma$.
\end{remark}

\begin{PROOF}{\ref{d43}}
1) Recalling $\lambda \le 2^\chi$ and $\chi^{< \theta_\iota} +
(\sigma_2)^{< \theta_\iota} \le \sigma_1$ hence 
$\chi^{< \theta_\iota} + 2^{< \theta_\iota} \le \sigma_1$, choose
\mn
\begin{enumerate}
\item[$(*)_1$]  
\begin{enumerate}
\item[$(a)$]  $A_{\alpha} \subseteq \chi$ 
(for $\alpha < \lambda$) which are pairwise distinct.
\sn
\item[$(b)$]  Let $\{(a_{i},d_{i}):i < \sigma_1\}$ be a list (maybe with
  repetitions) of the pairs $(a,d)$ 
satisfying $a \subseteq \chi,|a| < \theta_\iota$ and 
$d$ a function from $\cP(a)$ to $\sigma_2$ such that

\[
|\{b:b \subseteq a \text{ and } d(b) \ne 0\}| < \theta_\iota.
\]
\end{enumerate}
\end{enumerate}
\mn
Choose
\mn
\begin{enumerate}
\item[$(*)_2$]  $\bold c$ to be a symmetric two-place function 
from $\lambda$ to $\sigma_1$ exemplifying

\[
\Per_1(\lambda,\mu,\sigma_1,\bar\theta).
\]
\end{enumerate}
\mn
Now we define the two place function $\bold d$ from $\lambda$ to
$\sigma_2$ as follows: for $\alpha_0 < \alpha_1$:

\[
\bold d(\alpha_0,\alpha_1) = \bold d(\alpha_1,\alpha_0) := 
d_{\bold c(\alpha_0,\alpha_1)}
(A_{\alpha_\iota} \cap a_{\bold c(\alpha_0,\alpha_1)}).
\]

\mn
We shall show that $\bold d$ witnesses
$\Qr_{0,\iota}(\lambda,\mu,\sigma_2,\bar\theta)$ thus finishing upon
using Observation \ref{d40}(1) which yields the parallel assertion
about $\Pr_{0,\iota}(\lambda,\mu,\sigma_2,\bar\theta)$ 
because its assumption on the cardinals follows from those
of \ref{d43}(1), i.e. recall $\lambda = \mu = \cf(\mu)$ and $\theta_0
+ \theta_1 < \lambda$ so $\theta_\iota < \cf(\mu)$ and 
$\sigma^{< \theta_\iota}_2 < \mu$.
So let $\langle t_{\alpha}:\alpha<\mu\rangle$ be pairwise disjoint subsets of 
$\lambda,t_\alpha = t^0_\alpha \cup t^1_\alpha$ and 
$h^\iota_\alpha:t^\iota_\alpha \rightarrow \sigma_2$ such that
$h^{1-\iota}_\alpha$ is constant, $|t^0_\alpha| <
\theta_0,|t^1_\alpha| < \theta_1$ and $\pr:\sigma_2 \times \sigma_2
\rightarrow \sigma_2$.   As $\lambda = \mu = \cf(\mu)$
\wilog \, $\alpha < \beta < \mu \Rightarrow \sup(t_\alpha) < \min(t_\beta)$.
We have to find $\alpha_0 < \alpha_1$ as in the definition of
$\Qr_{0,\iota}(\lambda,\mu,\sigma_\iota,\bar\theta)$ 
see Definition \ref{d35}.  
As by assumption $\mu = \cf(\mu) > \theta$ and, of course, $\alpha <
\mu \wedge \ell < 2 
\Rightarrow \otp(t^\ell_\alpha) < \theta_\ell \le \theta$ \wilog \,  
there are $\varp^*_0 < \theta_0,\varp^*_1 < \theta_1$ such that
$\bigwedge\limits_{\alpha} \otp(t^\ell_\alpha) =
\varepsilon^*_\ell$ for $\ell = 0,1$.

For each $\alpha < \mu$ and $\ell < 2$ let $t^\ell_{\alpha} =
\{\zeta^\ell_{\alpha,\varepsilon}:\varepsilon < \varepsilon^*_\ell\}$ with
$\zeta^\ell_{\alpha,\varepsilon}$ increasing with $\varepsilon$.  As
$|\{\langle h^\iota_\alpha(\zeta^\iota_{\alpha,\varepsilon}):\varepsilon <
\varepsilon^*_\iota\rangle:\alpha < \mu\}| 
\le \sigma^{|\varepsilon^*_\iota|}_2 \le \sigma_2^{< \theta_\iota}  
< \mu = \cf(\mu)$, \wilog \, 
$h^\iota_\alpha(\zeta^\iota_{\alpha,\varepsilon}) = \xi^\iota_\varepsilon <
\sigma_2$ for all $\varepsilon < \varepsilon^*_\iota$ and
$h^{1-\iota}_\alpha(\zeta^{1-\iota}_{\alpha,\varp}) =
\xi^{1-\iota}_\varp$ which does not depend on $\alpha$.  
Renaming \wilog \, $\pr(\xi^0_{\varp(0)},\xi^1_{\varp(1)}) =
\xi_{\varp(\iota)}$, so rename it $\xi_{\varp(\iota)}$ for $\varp(0) <
\varp^*_0,\varp(1) < \varp^*_1$.

We should find $\alpha_0 < \alpha_1 < \mu$ such that for
$\varepsilon_0 < \varepsilon^*_0,\varepsilon_1 < \varepsilon^*_1$ we have
$\zeta_{\alpha_0,\varepsilon_0} < \zeta_{\alpha_1,\varepsilon_1}$
(which follows) and $\bold d(\zeta^0_{\alpha_0,\varepsilon_0},
\zeta^1_{\alpha_1,\varepsilon_1}) =
\pr(h^0_{\alpha_0}(\zeta^0_{\alpha_0,\varp_0}),
h^1_{\alpha_1}(\zeta^1_{\alpha_1,\varp_1}))$
which is equal to $\pr(\xi^0_{\varp_0},\xi^1_{\varp_1})$. 
Choose $a_{\alpha} \subseteq \chi,|a_{\alpha}| =
|\varepsilon^*_\iota| < \theta_\iota$ such that
$\langle A_{\zeta^\iota_{\alpha,\varepsilon}} \cap a_\alpha:
\varepsilon < \varepsilon^*_\iota \rangle$ is a sequence of pairwise
distinct subsets of $a_\alpha$.
As $\cf(\mu) = \mu > \chi^{<\theta_\iota}$ without loss of generality for every
$\alpha < \lambda = \mu$ we have 
$a_{\alpha} = a^*$ and $A_{\zeta^\iota_{\alpha,\varepsilon}}
\cap a^* = a^*_\varepsilon$ for all $\varepsilon < \varepsilon^*_\iota$.

For some $i < \sigma_1$ we have $a_i = a^*$ and $d_{i}(a^*_\varepsilon) = 
\xi_\varepsilon$ for every $\varepsilon < \varepsilon^*_\iota$.
By the choice of $\bold c$ for some $\alpha_0 < \alpha_1 < \mu$ the function
$\bold c \rest t_{\alpha_0} \times t_{\alpha_1}$ is constantly $i$,
so $\varepsilon_0 < \varepsilon^*_0 \wedge \varepsilon_1 <
\varepsilon^*_1 \Rightarrow \bold c(\zeta^0_{\alpha_0,\varepsilon_0},
\zeta^1_{\alpha_1,\varepsilon_1})= i$, hence for every
$(\varepsilon_0,\varepsilon_1) \in \varepsilon^*_0 \times
\varepsilon^*_1$ we have

\[
\bold d(\zeta^0_{\alpha_0,\varepsilon_0},\zeta^1_{\alpha_1,\varepsilon_1}) =
d_i(A_{\zeta^\iota_{\alpha_\iota,\varp_\iota}} \cap a_{i}) = 
d_{i}(a^*_{\varepsilon_\iota}) = 
\xi_{\varepsilon_\iota} = \pr(h^0_{\alpha_0}
(\zeta^0_{\alpha_0,\varp_0}),h^1_{\alpha_1}(\zeta^1_{\alpha_1,\varp_1})) 
\]

\mn
as required.

\noindent
1A) Similarly.

\noindent
2) Similar to part (3), see remarks inside its proof.

\noindent
3) Let $\theta = \theta_0 + \theta_1$ but for part (2) we let $\theta
= \theta_\ell$ and let $\bold c_1:[\lambda]^2 \rightarrow \sigma$ witness
$\Pr_0(\lambda,\lambda,\sigma,\bar\theta)$ and let $\bar f = \langle
f_\alpha:\alpha < \lambda\rangle$ be such that $f_\alpha$ is a
one-to-one function from $\sigma$ onto $\sigma + \alpha$. 
Let $\langle A_\alpha:\alpha < \lambda\rangle$ be a sequence of
pairwise distinct subsets of $\sigma$ and let $\langle (a_i,d_i):i <
\sigma\rangle$ list the pairs $(a,d)$ such that $a \in 
[\sigma]^{< \theta},d:\cP(a) \times \cP(a) \rightarrow \sigma$ and 
$|\{(b_1,b_2):b_1 \subseteq a,b_2 \subseteq a$ and $\bold
c_1(b_1,b_2) \ne 0\}| < \theta$; for part (2) we use 
$d:\cP(a) \rightarrow \sigma$.

Now we define
$\bold c_2:[\lambda]^2 \rightarrow \lambda$ as follows: for $\alpha <
\beta < \lambda$ let $\bold c_2(\{\alpha,\beta\}) = 
f_\beta((d_{\bold c_1(\{\alpha,\beta\})}(A_\alpha \cap a_{\bold
  c_1(\{\alpha,\beta\})},A_\beta \cap a_{\bold
  c_1(\{\alpha,\beta\})})))$.

So let $\bar\zeta^\iota = \langle \zeta^\iota_{\alpha,i}:\alpha <
\lambda,i < \bold i_\iota\rangle$ for $\iota < 2$ and $h:\bold i_0
\times \bold i_1 \rightarrow \lambda$ be as in Definition
\ref{d32}(1) but for part (2), $h:\bold i_\ell \rightarrow \lambda$,
see \ref{d32}(2).  For $\iota = 0,1$ for each $\alpha < \lambda$ and
$i < \bold i_\iota$ we can
find $a_{\alpha,\iota} \in [\sigma]^{< \theta_\iota}$ such that
$\bar b_{\alpha,\iota} := \langle A_{\zeta^\iota_{\alpha,i}} \cap
a_{\alpha,\iota}:i < \bold i_\iota\rangle$ is a sequence of pairwise
distinct sets.

\Wilog \, $\alpha < \lambda \wedge \iota < 2 \Rightarrow
a_{\alpha,\iota} = a_\iota,\bar b^\iota_\alpha = \bar b_\iota$; also \wilog
\, $\sup(\Rang(h)) \subseteq \min\{\zeta^\iota_{\alpha,i}:\alpha <
\lambda,i < \bold i_\iota$ and $\iota < 2\}$.

Next let $\bar\beta^\iota_\alpha = \langle \beta^\iota_{\alpha,i_0,i_1}:i_0
< \bold i_0$ and $i_1 < \bold i_1\rangle$ be a sequence of ordinals $<
\sigma$ such that $f_{\zeta^\iota_{\alpha,i_1}}(\beta^\iota_{\alpha,i_0,i_1}) =
h(i_0,i_1)$ and \wilog \, $\bar\beta^\iota_\alpha = \bar\beta^\iota$;
actually for part (3) we use only $f_{\zeta^\iota_{\alpha,i_1}}$ 
but for part (2) we use $f_{\zeta^\iota_{\alpha,i_\iota}}$ for
the $\iota$ from there.

Let $a = a_0 \cup a_1$ so $a \in [\sigma]^{<(\theta_0 + \theta_1)}$
and let $d:\cP(a) \times \cP(a) \rightarrow \sigma$ be such that
$d(b^0_{i_0},b^1_{i_1}) = \beta^1_{i_0,i_1}$ and $d(b_0,b_1) = 0$ if
$b_0,b_1 \subseteq a$ and $(b_0,b_1) \notin
\{(b^0_{i_0},b^1_{i_1}):i_0 < \bfi_0,i_1 < \bfi_1\}$.  Let $j < \sigma$ be
such that $(a_j,d_j) = (a,d)$.

Lastly, by the choice of $\bold c_1$
we can find $\alpha < \beta$ such that $i_0 < \bold i_0 \wedge
i_1 < \bold i_1 \Rightarrow \bold
c_1(\{\zeta^0_{\alpha,i_0},\zeta^1_{\alpha,i_1}\}) = j$; and now
check.

\noindent
4) Similarly to the proof of part (3).

\noindent
5) As $\Pr_1(\lambda,\lambda,\lambda,\bar\theta)$ by monotonicity we
have $\Pr_1(\lambda,\lambda,\partial,\bar\theta)$ hence by part (4) we
have $\Pr_0(\lambda,\lambda,\partial,\bar\theta)$ and now by part (3)
we can deduce $\Pr_0(\lambda,\lambda,\lambda,\bar\theta)$ as promised.
\end{PROOF}
\bigskip

\centerline {$* \qquad * \qquad *$}
\bigskip

In Juhasz-Shelah \cite{JuSh:1025} we use $\Col(\lambda,\kappa)$,
i.e. $\Pr_0(\lambda,\lambda,2,\kappa^+)$ quoting \cite[Ch.III,\S4]{Sh:g}
that e.g. $(\lambda,\kappa) = ((2^{\aleph_0})^{++} + \aleph_4,\aleph_0)$ is
O.K.  But in fact less suffices (see Definition \ref{d32}).
\begin{conclusion}  
\label{d47}
1) For $\lambda = \kappa^{+4}$ we have
$\Pr_{1}(\lambda,\lambda,\lambda,\kappa^+)$ which implies
$\Pr_{0,0}(\lambda,\lambda,\lambda,(\aleph_0,\kappa^+))$ and hence trivially
$\Pr_{0,0}(\lambda,\lambda,2,(\aleph_0,\kappa^+))$ holds.

\noindent
2) If $\Pr_{0,0}(\lambda,\lambda,2,(\aleph_0,\kappa^+))$ or just
$\Pr^{\uf}_{0,0}(\lambda,\lambda,2,(\aleph_0,\kappa^+))$,
 e.g. $\lambda = \aleph_4,\kappa = \aleph_0$  \then \, we have:
\mn
\begin{enumerate}
\item[$(*)_{\lambda,\kappa}$]  there is a topological space $X$ such
  that
\sn
\begin{enumerate}
\item[$(a)$]  $X$ is $T_3$, even has a clopen basis and has weight $\le
  \lambda$
\sn
\item[$(b)$]  the closure of any set of $\le \kappa$ points is compact
\sn
\item[$(c)$]  any infinite discrete set has an accumulation point
\sn
\item[$(d)$]  the space is not compact
\sn
\item[$(e)$]  some non-isolated point is not the accumulation point of
  any discrete set.
\end{enumerate}
\end{enumerate}
\end{conclusion}

\begin{PROOF}{\ref{d47}}
1) First we apply Theorem \ref{e8} (or \cite[Ch.III,\S4]{Sh:g}) with 
$(\kappa^{+4},\kappa^{+3},\kappa^+)$ here standing for 
$(\lambda,\partial,\theta)$ there.  Clearly the assumptions there hold
hence $\Pr_1(\kappa^{+4},\kappa^{+4},\kappa^{+4},\kappa^+)$ holds.

Second, we apply Claim \ref{d43}(1) with $0,\kappa^{+4},
\kappa^{+4},\kappa^{+3},\kappa^{+3},\kappa^+,\aleph_0,\kappa^+,\kappa^{+3}$
here standing for $\iota,\lambda,\mu,\sigma_1,\sigma_2,
\theta,\theta_0,\theta_1,\chi$ there.
Clearly the assumptions there hold because:
\mn
\begin{enumerate}
\item[$\bullet_1$]  ``$\Pr_1(\lambda,\mu,\sigma_1,\bar\theta)"$ there means
$\Pr_1(\kappa^{+4},\kappa^{+4},\kappa^{+3},(\aleph_0,\kappa^+))$ 
here which holds by the ``first" above and monotonicity
\sn
\item[$\bullet_2$]    ``$\chi^{< \theta_\iota} < \mu"$ there 
means ``$(\kappa^{+3})^{< \aleph_0} < \kappa^{+4}"$
\sn
\item[$\bullet_3$]   ``$\chi^{< \theta_\iota} \le \sigma_1"$ there means
``$(\kappa^{+3})^{< \aleph_0} \le \kappa^{+3}"$
\sn
\item[$\bullet_4$]   ``$2^\chi \ge \lambda"$ there means 
``$2^{\kappa^{+3}} \ge \kappa^{+4}$"
\sn
\item[$\bullet_5$]  $``\sigma^{< \theta_\iota}_2 \le \sigma_1"$ 
there which means here ``$(\kappa^{+3})^{< \aleph_0} \le \kappa^{+3}"$
\sn
\item[$\bullet_6$]  ``$\sigma^{< \theta_\iota}_2 < \mu"$ there
which means here $(\kappa^{+3})^{< \aleph_0} < \kappa^{+4}$
\end{enumerate}
\mn
So all of them hold indeed.

Next, the conclusion of \ref{d43}(1) is
$\Pr_{0,\iota}(\lambda,\mu,\sigma_2,\bar\theta)$ which here 
means $\Pr_{0,0}(\kappa^{+4},\kappa^{+4},\kappa^{+3},(\aleph_0,\kappa^+))$.

Lastly, by \ref{d43}(2) we get
$\Pr_{0,0}(\kappa^{+4},\kappa^{+4},\kappa^{+4},(\aleph_0,\kappa^+))$.

\noindent
2) By Claim \ref{d53} below, which generalize the proof of 
Juhasz-Shelah \cite{JuSh:1025}, that is, 
let $\bar D = \langle D_i:i < \beth_2\rangle$ list the
ultrafilters on $\sigma := \aleph_0$ and let $\sigma_i = \sigma$ for
$i < \beth_2$ and $\theta = \kappa^+$.  So clause (A) of \ref{d53} below
holds, hence we can apply \ref{d53} for 
$(\lambda,\theta) = (\lambda,\kappa^+)$ and $\bar D$.  So clause (a)
of \ref{d47}(2) holds by $(B)(a)(\alpha)$ of \ref{d53}, of course; clause
(b) of \ref{d47}(2) holds by $(B)(a)(\gamma)$ recalling the choice of
$\bar D$; clause (c) there holds by $(B)(a)(\varp)$; clause
(d) there holds by $(B)(a)(\delta)$; and lastly, clause (e) there 
holds by $(B)(b)$.  So we are done.
\end{PROOF}

\noindent
Moreover
\begin{claim}
\label{d50}
1) If $\kappa$ is regular and $\lambda = \kappa^{+3}$ then
$\Pr_{1}(\lambda,\lambda,\lambda,(\aleph_0,\kappa^+))$ hence
$\Pr_{0,0}(\lambda,\lambda,\lambda,(\aleph_0,\kappa^+))$.

\noindent
2) $(*)_{\aleph_3,\aleph_0}$ from \ref{d47}(2) holds.

\noindent
3) $(*)_{\kappa^{+3},\kappa}$ from \ref{d47}(2) holds for $\kappa$ regular.
\end{claim}

\begin{PROOF}{\ref{d50}}
Like the proof of \ref{d47} using Theorem \ref{e16} instead of Theorem
\ref{e8}, that is, we apply \ref{e16} with
$(\aleph_3,\aleph_2,\aleph_1,\aleph_0)$ standing for
$(\lambda,\partial,\theta_1,\theta_0)$.  
\end{PROOF}

\noindent
We conclude this section with an explicit proof of the topological
statement in \ref{d47}(2).  We shall need the following:
\begin{definition}
\label{d531}
Let $X$ be a topological space, $D$ an ultrafilter over $\sigma$.

\noindent
1) An element $y \in X$ is the $D$-limit of a sequence of points
$\langle x_j:j < \sigma\rangle$ in $X$ \Iff \, $y \in u \Rightarrow
\{j < \sigma:x_j \in u\} \in D$ whenever $u$ is a open subset of
$X$.

\noindent
2) $X$ is $D$-complete \Iff \, for every sequence of points
$\langle x_j:j < \sigma\rangle$ in $X$ there is $y \in X$ such that
$y$ is the $D$-limit of the sequence.

\noindent
3) If $\bar D = \langle D_i:i < i_*\rangle$ is a sequence such that
each $D_i$ is an ultrafilter over $\sigma_i = \sigma(i)$ \then \, $X$ is $\bar
D$-complete iff $X$ is $D_i$-complete for every $i < i_*$.
\end{definition}

\begin{claim}
\label{d53}
If (A) then (B) where
\mn
\begin{enumerate}
\item[(A)]
\begin{enumerate}
\item[(a)]  $\lambda = \cf(\lambda) > \theta = \cf(\theta) > \aleph_0$
\sn
\item[(b)]  $\bar D = \langle D_i:i < i_*\rangle$, each $D_i$ 
is a non-principal ultrafilter on $\sigma_i$ and $\sigma_i < \theta$
\sn
\item[(c)]  $\Pr_{0,0}(\lambda,\lambda,2,(\aleph_0,\theta))$; yes!
  $\Pr_{0,0}$ and not $\Pr_0$
\end{enumerate}
\sn
\item[(B)]  there is a topological space $X$ and a
 point $g \in X$ such that:
\sn
\begin{enumerate}
\item[(a)]
\begin{enumerate}
\item[$(\alpha)$]  $X$ is a subspace of ${}^\lambda 2$ hence has a
  clopen basis and is a $T_3$-space
\sn
\item[$(\beta)$]  $X$ is a dense subset of ${}^\lambda 2$ hence has no
  isolated point and its weight is $\lambda$
\sn
\item[$(\gamma)$]  if every non-principal ultrafilter $D$ on
  a cardinal $\sigma < \theta$ appears in $\bar D$ \then \, for any set $Y
  \subseteq X$ of cardinality $< \theta$, the closure of $Y$ is
  compact
\sn
\item[$(\delta)$]  $X$ is not compact
\sn
\item[$(\varp)$]  any subset of $X$ of cardinality $\ge
  \min\{\sigma_i:i < i_*\}$ has an accumulation point; so the
  cardinality can be $\aleph_0$
\sn
\item[$(\zeta)$]  $X$ is $\bar D$-complete
\end{enumerate}
\sn
\item[(b)]
\begin{enumerate}
\item[$(\alpha)$]  $g \in X$ is not an accumulation point of any
  discrete set $Y \subseteq X \backslash \{g\}$
\sn
\item[$(\beta)$]  moreover, $g$ is not an accumulation point
of any set $Y \subseteq X \backslash \{g\}$ of cardinality $< \lambda$
\end{enumerate}
\sn
\item[(c)]  
\begin{enumerate}
\item[$(\alpha)$]  $X$ has $\le \lambda^{< \theta} +
  \sum\limits_{\sigma < \theta} 2^{2^\sigma}$ points
\sn
\item[$(\beta)$]  $X$ has $\ge \lambda$ points
\end{enumerate}
\sn
\item[(d)]  if $i_* < \lambda$ and $\alpha < \lambda \Rightarrow
  |\alpha|^{< \theta} < \lambda$ \then \,
\sn
\begin{enumerate}
\item[$(\alpha)$]  $X$ has no discrete subset of cardinality $\ge
  \lambda$, moreover
\sn
\item[$(\beta)$]  $h L^+(X) \le \lambda$ so $\lambda = \mu^+
  \Rightarrow h L(X) \le \mu$.
\end{enumerate}
\end{enumerate}
\end{enumerate}
\end{claim}

\begin{PROOF}{\ref{d53}}
\medskip

\noindent
\underline{Stage A}:   We make some choices:
\mn
\begin{enumerate}
\item[$(*)_1$]  
\begin{enumerate}
\item[(a)]  let $\bold c:[\lambda]^2 \rightarrow \{0,1\}$ witness
  $\Pr_{0,0}(\lambda,\lambda,2,(\aleph_0,\theta))$
\sn
\item[(b)]  let $\bar h^* = \langle h^*_\alpha:\alpha <
  \lambda\rangle$ list the finite partial functions from $\lambda$ to
$\{0,1\}$; \wilog \, $\dom(h^*_\alpha) \subseteq \alpha$
\sn
\item[(c)]  let $g \in {}^\lambda 2$ be constantly 1.
\end{enumerate}
\end{enumerate}
\mn
Further
\mn
\begin{enumerate}
\item[$(*)_2$]  for $\alpha < \lambda$ we define $f^*_\alpha \in
  {}^\lambda 2$ as follows:
\sn
\begin{itemize}
\item  for $\beta < \lambda$ we let $f^*_\alpha(\beta)$ be
\begin{enumerate}
\item[(a)]  $h^*_\alpha(\beta)$ if $\beta \in \dom(h^*_\alpha)$
\sn
\item[(b)]  $\bold c\{\beta,\alpha\}$ if $\beta < \alpha \wedge \beta
  \notin \dom(h^*_\alpha)$
\sn
\item[(c)]  0 otherwise, i.e. if $\beta \ge \alpha$.
\end{enumerate}
\end{itemize}
\end{enumerate}
\mn
Our $X$ will include each $f^*_\alpha$ for $\alpha < \lambda$ but more.
\mn
\begin{enumerate}
\item[$(*)_3$]  for $\beta \le \lambda$ we let
\sn
\begin{enumerate}
\item[(a)]  $\cF_\beta = \{f^*_\alpha:\alpha < \beta\}$
\sn
\item[(b)]  $\cF^*_\beta = c \ell_{\bar D}(\cF_\beta)$,
  i.e. $\cF^*_\beta$ is the minimal subset of ${}^\lambda 2$ which
  includes $\cF_\beta$ and is $\bar D$-closed
\sn
\item[(c)]  $\cG^*_\beta = \{f:f \in \cF^*_\lambda$ and $f \rest
  [\beta,\lambda)$ is constantly zero$\}$.
\end{enumerate}
\end{enumerate}
\mn
So
\mn
\begin{enumerate}
\item[$(*)_4$]  $\cF^*_\lambda$ is the union of the
  $\subseteq$-increasing sequence $\langle \cF^*_\beta:\beta <
  \lambda\rangle$.
\end{enumerate}
\mn
[Why?  Clearly $\langle \cF_\beta:\beta < \lambda\rangle$ is
$\subseteq$-increasing and as $\cf(\lambda) \ge \theta$ and $D_i$ is
an ultrafilter on $\sigma_i < \theta$ for $i < i_*$ clearly $(*)_4$ follows.]

Lastly, we choose $X$
\mn
\begin{enumerate}
\item[$(*)_5$]  $X$ is the subspace of ${}^\lambda 2$ with set of
  elements $\cF^*_\lambda \cup \{g\}$.
\end{enumerate}
\mn
So it suffices to prove that $X,g$ are as required in the claim.
\mn
\begin{enumerate}
\item[$(*)_6$]  if $f \in \cF^*_\lambda$ \then \, for some triple
  $(u,v,D)$ we have: 
\sn
\begin{enumerate}
\item[(a)]  $u,v \in [\lambda]^{< \theta}$
\sn
\item[(b)]  $D$ an ultrafilter on $u$
\sn
\item[(c)]  $f = \lim_D(\langle f^*_\alpha:\alpha \in u\rangle)$
\sn
\item[(d)]  if $\beta \in \lambda \backslash v$, \then \, $f(\beta) =
  1 \Leftrightarrow \{\alpha \in u:\beta < \alpha$ and $\bold
  c\{\alpha,\beta\}=1\} \in D$. 
\end{enumerate}
\end{enumerate}
\mn
[Why?  Recall $\cF^*_\lambda$ is $c \ell_{\bar D}(\cF_\lambda)$ and each
$D_i$ is an ultrafilter on some $\sigma < \theta$.  Hence we can find
a sequence $\langle f^*_\alpha:\alpha \in [\lambda,\alpha_*)\rangle$
listing $\cF^*_\lambda \backslash \cF_\lambda$ and for each such
$\alpha,i(\alpha) = i_\alpha < i_*$ and $\bar\beta_\alpha \in
{}^{\sigma(i(\alpha))} \lambda$ are such that $f^*_\alpha =
\lim_{D_{i(\alpha)}}(\langle f_{\beta_\alpha,\varp}:\varp <
\sigma_{i(\alpha)}\rangle)$.  As $\theta$ is regular, \underline{clearly} 
there are $u \in [\lambda]^{< \theta}$ and an ultrafilter $D$ 
on $u$ such that clause (c) holds.

Why?  If $f = f^*_\alpha,\alpha < \lambda$ then $u = \{\alpha\}$ is
as required and if $f=f^*_\alpha,\alpha \in [\lambda,\alpha_*)$ then we can
prove this by induction on $\alpha$.

Now choose $v = \cup\{\dom(h^*_\alpha):\alpha \in u\}$, clearly $u,v$
are as required.  E.g. if $f = f^*_\alpha,\alpha < \lambda$ the
ultrafilter $D$ is the unique principal ultrafilter on $\{\alpha\}$;
for $(*)_6(d)$ recall the choice of the $f^*_\alpha$'s for $\alpha < \lambda$.]
\mn
\begin{enumerate}
\item[$(*)_7$]  if $f \in \cF^*_\lambda$ and $\delta < \lambda$ has
  cofinality $\ge \theta$, \then \, for some $\gamma < \delta$, at
  least one of the following holds: 
\sn
\begin{enumerate}
\item[(a)]  if $\beta \in [\gamma,\lambda)$ then $f(\beta) = 0$
\sn
\item[(b)]  for some $u = u_f \in [\lambda \backslash \delta]^{<
    \theta}$ and $v=v_f \in [\lambda \backslash \delta]^{< \theta}$
 and ultrafilter $D$ on $u$ we have
\sn
\begin{itemize}
\item  if $\beta \in [\gamma,\lambda) \backslash v_f$ 
then $f(\beta) = \lim_D(\langle \bold c\{\beta,\alpha\}:\alpha \in u\rangle)$.
\end{itemize}
\end{enumerate}
\end{enumerate}
\mn
[Why?  Let $u,v,D$ be as in $(*)_6$.  If $u \cap \delta \in D$ then
let $\gamma$ be $\sup(u \cap \delta) < \delta$ and by
$(*)_2(c) + (*)_6(c)$ clearly clause (a) of $(*)_7$ holds.  So we can
assume $u \cap \delta \notin D$ and as $D$ is an ultrafilter on $u$,
necessarily $u \backslash \delta \in D$.  Let $u' = u \backslash
\delta,\gamma = \sup(\cup\{\dom(h^*_\alpha) \cap \delta:\alpha \in u\}
\cup (v \cap \delta)) +1$ and $D' = D \cap \cP(u')$ and $v' = v
\backslash \delta$, they clearly witness clause (b) of
$(*)_7$.  Together we are done.]
\mn
\begin{enumerate}
\item[$(*)_8$]
\begin{enumerate}
\item[(a)]   if $f \in \cF^*_\lambda$, \then \, for some $\beta <
  \lambda$ we have $f \in \cF^*_\beta$ which implies $f$ is constantly
  zero on $[\beta,\lambda)$
\sn
\item[(b)]    $\cF^*_\beta \subseteq \cG^*_\beta \subseteq \cF^*_\lambda$
\sn
\item[(c)]   $\cG^*_\beta$ is $\subseteq$-increasing with $\beta$ with
  union $\cF^*_\lambda$.
\end{enumerate}
\end{enumerate}
\mn
[Why?  Clause (a) holds by $(*)_3(b) + (*)_4$ above.  Clauses (b),(c)
are easy too recalling $(*)_3(a)$.]
\medskip

\noindent
\underline{Stage B}:  Now we check the demands in (B) of the claim.
\mn
\begin{enumerate}
\item[$\oplus_1$]   $X$ is a subspace of ${}^\lambda 2$ [so clause
  $(B)(a)(\alpha)$ holds] hence $X$ is a $T_3$ topological space 
with a clopen base.
\end{enumerate}
\mn
[Why?  By its choice in $(*)_5$.]
\mn
\begin{enumerate}
\item[$\oplus_2$]   $X$ is dense in ${}^\lambda 2$ hence clause
  $(B)(a)(\beta)$ holds.
\end{enumerate}
\mn
[Why?  By the choice of $\bar h^*$ in $(*)_1(b)$ because $h^*_\alpha \subseteq
f^*_\alpha$ for $\alpha < \lambda$ by $(*)_2(a)$.]
\mn
\begin{enumerate}
\item[$\oplus_3$]   $X$ is $D_i$-complete for every $i < i_*$ hence
  clause $(B)(a)(\zeta)$ holds.
\end{enumerate}
\mn
[Why?  By the choice of $\cF^*_\lambda$ in $(*)_3(b)$ because $X \backslash
\cF^*_\lambda = \{g\}$ recalling $\lambda = \cf(\lambda) > \theta$.]
\mn
\begin{enumerate}
\item[$\oplus_4$]   $\lambda \le |X| \le \lambda^{< \theta} +
  \sum\limits_{\sigma < \theta} 2^{2^\sigma}$ and also $|X| \le
  \lambda^{< \theta} + 2^{\theta + |i_*|}$ hence clause (B)(c) holds.
\end{enumerate}
\mn
[Why?  Clearly $|\cF_\lambda|=\lambda$ and $\cF_\lambda \subseteq
\cF^*_\lambda \subseteq X$ hence $\lambda \le |X|$.  As $|X \backslash
\cF^*_\lambda| = |\{g\}| =1$ and by $(*)_6$ the other inequalities follow.]
\mn
\begin{enumerate}
\item[$\oplus_5$]   $g \notin c \ell(Y)$ \when \, $Y \subseteq X
  \backslash \{g\}$ and at least one of the following holds:
\begin{enumerate}
\item[(a)]  $|Y| < \lambda$
\sn
\item[(b)]  for some $\beta < \lambda,Y \subseteq \cF^*_\beta$
\sn
\item[(c)]  for some $\beta < \lambda,
Y \subseteq \cG^*_\beta := \{f \in \cF^*_\lambda:f \rest
[\beta,\lambda]$ is constantly zero$\}$.
\end{enumerate}
\end{enumerate}
\mn
[Why?  If clause (a), i.e. $|Y| < \lambda = \cf(\lambda)$ as
  $\langle\cF^*_\beta:\beta < \lambda\rangle$ is
  $\subseteq$-increasing with union $\cF^*_\lambda$ by $(*)_4$, necessarily $Y
  \subseteq \cF^*_\beta$ for some $\beta < \lambda$, i.e. clause (b); but
this in turn implies clause (c) by $(*)_8(b)$.

But if clause (c) holds for $\beta$, then $g \notin c \ell(Y)$
recalling that $g(\gamma) = 1$ for every $\gamma < \lambda$.]

Now comes a major point using the choice of $\bold c$,
i.e. $\Pr_{0,0}(\lambda,\lambda,2,(\aleph_0,\theta))$. 
\mn
\begin{enumerate}
\item[$\oplus_6$]   if $Y \subseteq \cF^*_\lambda$ and $\beta < \lambda
  \Rightarrow Y \nsubseteq \cG^*_\beta$ \then \, $Y$ is
 not discrete and even not left separated (hence, together with
 $\oplus_5$, clause (B)(b) holds).
\end{enumerate}
\mn
[Why?  For $\alpha < \lambda$ choose $f_\alpha \in Y \backslash
\cG^*_\alpha \subseteq \cF^*_\lambda \backslash \cF_\alpha$ hence 
there is $\beta^1_\alpha \in [\alpha,\lambda)$ such that 
$f_\alpha(\beta^1_\alpha) = 1$ and there is $\beta^2_\alpha \in
(\beta^1_\alpha,\lambda)$ such that $f_\alpha \rest
[\beta^2_\alpha,\lambda)$ is constantly zero.

Recall that ``$Y$ is left separated (in the space $X$)" means that there is
a well-ordering $<^*$ on $Y$ such that for every $x \in Y$ the set
$\{y \in Y:x <^* y\}$ is closed in the induced topology on $Y$.

Toward contradiction assume $Y$ is discrete or just left
separated.  Fix a well-ordering $<^*$ on $Y$ which witnesses this fact.
Clearly we can find $\cU_0 \in [\lambda]^\lambda$ such that
$\langle \beta^1_\alpha:\alpha \in \cU_0\rangle$ is an increasing
sequence of ordinals and on $Y,<^*$ and the usual order agree.

Now by the choice of $<^*$ for some $\cU \in [\cU_0]^\lambda$ 
we can find a sequence
$\bar h = \langle h_\alpha:\alpha \in
\cU\rangle,h_\alpha$ is a finite function from $\lambda$ to $\{0,1\}$
satisfying (the statement $\bullet_0 + \bullet_2$ by the definition of
``$<^*$ witness $Y$ is left separated"; the statements $\bullet_1$
hold as \wilog \, as increasing $h_\alpha$ makes no harm, and the
statement $\bullet_3$ holds \wilog \, because we can replace $\cU$ by
any $\cU' \in [\cU]^\lambda$):
\mn
\begin{enumerate} 
\item[$\bullet_0$]  $h_\alpha \subseteq f_\alpha$
\sn
\item[$\bullet_1$]   $\beta^1_\alpha,\beta^2_\alpha \in \Dom(h_\alpha)$
\sn
\item[$\bullet_2$]  if $\alpha_1 < \alpha_2$ then $h_{\alpha_1}
  \nsubseteq f_{\alpha_2}$.  Also (not used)
\sn
\item[$\bullet_3$]  if $\alpha_1 < \alpha_2$ are from $\cU$ then 
$\beta^2_{\alpha_1} < \beta^1_{\alpha_2}$ hence $h_{\alpha_2}
\nsubseteq f_{\alpha_1}$. 
\end{enumerate}
\mn
Renaming \wilog \,
\mn
\begin{enumerate}
\item[$\bullet_4$]  $\cU =\lambda$ and still $\beta^2_\alpha >
\beta^1_\alpha \ge \alpha,f_\alpha(\beta^1_\alpha) = 1$ and $f_\alpha
\rest [\beta^2_\alpha,\lambda)$ is constantly zero.
\end{enumerate}
\mn
For each $\delta \in S_1 := S^\lambda_\theta =
\{\delta < \lambda:\cf(\delta) = \theta\}$ we consider $(*)_7$ with
$(f_\delta,\delta)$ here standing for $(f,\delta)$ there, now
$\beta^1_\delta \ge \delta,f_\delta(\beta^1_\delta) = 1$ by
$\bullet_4$ hence clause $(*)_7(a)$ fails, so necessarily clause
$(*)_7(b)$ holds.
So there is a quadruple $(\gamma_\delta,u_\delta,v_\delta,D_\delta)$ as 
there\footnote{They depend also on $f=f_\delta$, but $\delta$
  determines $f$.}  and let $\beta^3_\delta :=
\sup(\delta \cap (\dom(h_\delta)))$, as $h_\delta$ is a finite function,
necessarily $\beta^3_\delta < \delta$.  So by Fodor lemma for
some $\gamma_* < \lambda$ the set $S_2 = \{\delta \in
S_1:\gamma_\delta,\beta^3_\delta \le \gamma_* < \delta\}$ 
is stationary hence so is $S_3 =
\{\delta \in S_2$: if $\alpha < \delta$ then $u_\alpha,v_\alpha \subseteq
\delta,\beta^1_\alpha < \delta,\beta^2_\alpha < \delta$ 
and $\dom(h_\alpha) \subseteq \delta\}$.  As $\dom(h_\alpha)$
is finite and $\range(h_\alpha) \subseteq \{0,1\}$ clearly for 
some $h_*,h_{**}$ the set $S_4 = \{\delta \in S_3:h_\delta
\rest \delta = h_*$ and $h_{**} = \{(\otp(\dom(h_\delta) \cap
\gamma),h_\delta(\gamma)):\gamma \in \dom(h_\delta)\}\}$ is stationary.

For $\delta \in S_4$ let 
$u_{\delta,0} = \Dom(h_\delta) \backslash \Dom(h_*),h'_\delta = h_\delta
\rest u_{\delta,0}$ and $u_{\delta,1} = u_\delta$ and recall $u_\delta
\cap \delta = \emptyset = v_\delta \cap \delta$, 
see $(*)_7(b)$.  Note that
$\Qr_{0,0}(\lambda,\lambda,2,(\aleph_0,\theta))$ holds, see Definition
\ref{d35}(1),(2) for $\iota=0$, now it holds because we are assuming
$\Pr_{0,0}(\lambda,\lambda,2,(\aleph_0,\theta))$ by \ref{d40}(1).  So
we can apply the definition of 
$\Qr_{0,0}(\lambda,\lambda,2,(\aleph_0,\theta))$ and the choice 
of $\bold c$ to $\langle (u_{\delta,0},u_{\delta,1}:\delta \in
S_4\rangle$ and $\langle h'_\delta:\delta \in S_4\rangle$.  So there
are $\delta_1,\delta_2$ such that:
\mn
\begin{enumerate}
\item[$\bullet_5$]    $\delta_1 < \delta_2$ are from $S_4$
\sn
\item[$\bullet_6$]   if $\alpha \in u_{\delta_1,0}$ 
and $\beta \in u_{\delta_2,1}$ then
 $\bold c\{\alpha,\beta\} = h'_{\delta_1}(\alpha)$.
\end{enumerate}
\mn
Next
\mn
\begin{enumerate}
\item[$\bullet_7$]   if $\alpha \in u_{\delta_1,0}$ then
$f_{\delta_2}(\alpha) = \lim_{D_{\delta_2}}(\langle 
\bold c \{\alpha,\beta\}:\beta
  \in u_{\delta_2,1} = u_{\delta_2}\rangle)$.
\end{enumerate}
\mn
[Why?  By the choice of
$(\gamma_{\delta_2},u_{\delta_2},D_{\delta_2},h_*,h_{**})$ that is
recalling $(*)_7(b)$ because $\alpha \in
u_{\delta_1,0} \Rightarrow \alpha \in \dom(h'_{\delta_1}) \Rightarrow
\alpha \ge \delta_1 \Rightarrow \alpha \ge \gamma_* \ge
\gamma_{\delta_2}$ and $\alpha \in u_{\delta_1,0} \cup v_{\delta_1}
\Rightarrow \alpha < \delta_2$.] 
\mn
\begin{enumerate}
\item[$\bullet_8$]   if $\alpha \in \dom(h'_{\delta_1})$ then
  $f_{\delta_2}(\alpha) = h'_{\delta_2}(\alpha)$.
\end{enumerate}
\mn
[Why?  By $\bullet_7$ because $u_{\delta_1,0} = \dom(h'_{\delta_1})$
and $\bullet_6$.]
\mn
\begin{enumerate}
\item[$\bullet_9$]  $h'_{\delta_1} \subseteq f_{\delta_2}$.
\end{enumerate}
\mn
[Why?  By $\bullet_8$.]

However, $h_{\delta_1} \subseteq f_{\delta_1}$ by $\bullet_0$ hence
$h_* \subseteq h_{\delta_1} \subseteq f_{\delta_1}$ but $h_* \subseteq
h_{\delta_2} \nsubseteq f_{\delta_1}$ by $\bullet_2$ and
$h'_{\delta_2} = h_{\delta_2} \rest (\dom(h_{\delta_2}) \backslash
\dom(h_*)$ hence
\mn
\begin{enumerate}
\item[$\bullet_{10}$]   $h'_{\delta_2} \nsubseteq f_{\delta_1}$.
\end{enumerate}
\mn
But $\bullet_{10}$ contradict $\bullet_9$, all this follows from the
assumption toward contradiction in the beginning of the proof of
$\oplus_6$, so $\oplus_6$ holds indeed.

Now we can check all the remaining demands in (B), e.g.
\medskip

\noindent
\underline{Clause $(B)(d)(\beta)$}:  Assume toward contradiction that $h
L^+(X) > \lambda$.  This means that some $Y \subseteq X$ has
cardinality $\lambda$ and is right separated (by some well ordering).  Now
\wilog \, $g \notin Y$ and if $\beta < \lambda \Rightarrow Y
\nsubseteq \cG^*_\beta$ then we get a contradiction by $\oplus_6$.  So
we are left with the case $Y \subseteq \cG^*_\beta$ for some $\beta <
\lambda$.  But by the clause assumption 
$|\cG^*_\beta| \le |\beta|^{< \theta} + |i_*|$ which has
cardinality $< \lambda$, so we are done proving $(B)(d)(\beta)$. 

We are done proving \ref{d53}: most clauses of 
(B) were proved and we have to add that:
clauses $(B)(a)(\gamma) + (\varp)$ hold by the choice of
$\cF^*_\lambda$ as $X \backslash \cF^*_\lambda = \{g\}$.  Clause
$(B)(a)(\delta)$ is exemplified by any uniform ultrafilter $D$ on
$\lambda$ such that $\{\alpha:f^*_\alpha(0)=r\} \in D$, exists by
$(*)_3(c) + (*)_8$.
\end{PROOF} 
\newpage

\section {The colouring existence}

We try to explain the proof of \ref{e4}, \ref{e16}; probably more of
it will make sense after reading part of the proof.

Claim \ref{e4} should be understood as follows: given a set $S$
and functions $F_\iota:S \rightarrow \kappa_\iota$ for 
$\iota=0,1$ and a sequence $\varrho \in 
{}^{\omega >}S,\bold d(\varrho)$ is a natural number which in the
interesting case is a ``place in the sequence", i.e. $\bold d(\varrho)
< \ell g(\varrho)$.

In the interesting cases, $\varrho = \eta_0 \char 94 \nu_0 \char 94
\rho \char 94 \nu_1 \char 94 \eta_1$ is as constructed during the proof of
\ref{e16}, and if (B)(a)-(d) of \ref{e4} holds, $\ell g(\eta_0) +
\ell_4$ is a place in the sequence; so \ref{e4} tells us that it
depends only on $\varrho$ (and not on the representation
$(\eta_0,\nu_0,\rho,\nu_1,\eta_1)$ of $\varrho$).

How does $\bold d$ help us in the proof of Theorem \ref{e16}?

We shall describe it for the case of $\theta_1$ colours, i.e. $\sigma
= \theta_1$ and the colouring is called $\bold c_1$.  
Let $(\kappa_0,\kappa_1,\kappa_2) = (\theta_0,\theta_1,\lambda)$.
We shall be given pairwise disjoint
$t_\alpha = t^0_\alpha \cup t^1_\alpha$ for $\alpha
< \lambda$ and a colour $j_* < \theta_1$ such that $|t^\iota_\alpha| <
\theta_\iota$ for $\iota = 0,1$ and $\alpha < \lambda$ and 
we shall carefully choose
$\alpha_0 < \alpha_1$ exemplifying the desired conclusion.

Toward choosing the pair $(\alpha_0,\alpha_1)$ 
we also choose $\delta_0 < \delta_1 < \delta_2 < \delta_3$ which will be
from $(\alpha_0,\alpha_1)$ such that $\sup(t_{\alpha_0}) < \delta_0$
and $\ell_4$ such that:
\mn
\begin{enumerate}
\item[$(a)$]  we let $\nu_0 = \rho_{\bar h}(\delta_3,\delta_2),\rho =
  \rho_{\bar h}(\delta_2,\delta_1),\nu_1 = \rho_{\bar
  h}(\delta_1,\delta_0)$ where $\rho_{\bar h}(\delta',\delta'')$ is
  derived from the sequence $\rho(\delta',\delta'')$, see before
  $\odot_2$ in the proof of \ref{e16}
\sn
\item[$(b)$]   $\ell_4 < \ell g(\nu_0)$ and $h'(F_1(\nu_0(\ell_4))) =
  j_*$ where $h':\kappa_1 \rightarrow \kappa_2$ is chosen in $\odot_7$
  in the proof \ref{e16}
\sn
\item[$(c)$]   let $\zeta_0 \in t^0_{\alpha_0}$ and $\zeta_1 \in
  t^1_{\alpha_1}$ and define $\eta_{1,\zeta_0} = 
\rho_{\bar h}(\delta_0,\zeta_0),\eta_{0,\zeta_1} = \rho_{\bar
  h}(\zeta_1,\delta_3)$
\sn
\item[$(d)$]  continuing clause (c) by the construction
  $\varrho_{\zeta_1,\zeta_0} := \rho_{\bar h}(\zeta_1,\zeta_0)$ is
  equal to $\eta_{0,\zeta_1} \char 94 \nu_0 \char 94 \rho \char 94
  \nu_1 \char 94 \eta_{1,\zeta_0}$.
\end{enumerate}
\mn
So naturally we choose the colouring $\bold c_1$ such that $\bold
c_1(\alpha_0,\alpha_1) = h'(F_1(\varrho(\ell g(\eta_0) + \ell_4)))$ and
\ref{e4} tells us that assuming (a)-(d) this will be $j_*$.  Note it
is desirable that in \ref{e4}, the sequences $\eta_0,\eta_1$ in a
sense have little influence on the result, as they vary, i.e. we like
to get $j_*$ for every $\zeta_0 \in t^0_{\alpha_0},\zeta_1 \in
t^1_{\alpha_1}$.

Why do we demand in clause (b), $h_2(F_1(\nu_0(\ell_4))) = j_*$ and
not simply $F_1(\nu_0(\ell_4)) = j_*$ and similarly when defining
$\bold c_1$ in $\odot_7$ in the proof?  
Because we do not succeed to fully control
$F_1(\nu_0(\ell_4))$, but just to place it in some stationary 
$S \subseteq \theta_1$, however we can use $\theta_1$ pairwise disjoint
stationary set and $h_1$ tells us which one.

When we choose $\alpha_0 < \alpha_1$ (in stage C of the proof) we
first choose a pair $\delta_1 < \delta_2$ hence $\rho$ (in $\oplus_0$
of the proof), then we choose
an ordinal $\delta_0 < \delta_1$ hence $\nu_1$ (in $\oplus_{0.1}$ of
the proof) then $\varepsilon_* \in s_{\delta_2} \subseteq \kappa_1$
after $\oplus_{0.2}$ of the proof, (see below)
large enough.  Only then using
$\varepsilon_*$ we choose $\delta_3$ and then $\alpha_1$
 (also after $\oplus_{0.2}$) hence $\eta_{0,\zeta}$ for $\zeta \in
 t^1_{\alpha_1}$.  Lastly, we choose
$\alpha_0 < \delta_0$ hence $\eta_{1,\zeta_0}$ 
for $\zeta_0 \in t^0_{\alpha_0}$.  Of course,
those choices are under some restrictions.  More specifically, (in
stage B) though not determining any of
$\eta_{0,\zeta_0},\nu_0,\rho,\nu_1,\eta_{1,\zeta_1}$ we restrict them in
some ways.

Earlier, we first in $(*)_1$ choose
$\cU^{\up}_1,\alpha^*_1,\varepsilon^{\up}_{1,1},\varepsilon^{\up}_{1,0}$
with the intention that $\alpha_1 \in \cU^{\up}_1$ ``promising" that
if $\alpha_1 \in \cU^{\up}_1$ then $\Rang(F_1(\eta_0)) \subseteq
\varepsilon^{\up}_{1,1} < \kappa_1$, i.e. $\zeta_1 \in t^1_{\alpha_1}
\Rightarrow \Rang(F_1(\eta_{0,\zeta_1})) \subseteq
\varepsilon^{\up}_{1,1}$, similarly in the further steps
below.  Second we do not ``know" for which $\varepsilon < \kappa$ we
shall use $S^{\kappa_1}_{\kappa_0,\varepsilon} \subseteq \kappa_1$, so
we consider all of them, i.e. in $(*)_2$ we choose
$\cU^{\up}_{2,\varepsilon},g_{2,\varepsilon},\
\gamma^*_\varepsilon,\alpha^*_{2,\varepsilon}$ satisfying
$g_{2,\varepsilon}:\cU^{\up}_{2,\varepsilon} \rightarrow \cU^{\up}_1$
such that later $\delta_3 \in \cU^{\up}_{2,\varepsilon}$ and $\alpha_1
= g_{2,\varepsilon}(\delta)$.  We still do not know what $\nu_2$ will
be hence how to compute $\ell_4$, but $\rho_{\bar
  h}(\alpha_1,\delta_3)$ will be part of it and for each $\varepsilon
< \kappa_1$ we can compute $\ell_{2,\varepsilon}$ which will be the
first place $\ell$ in $\nu_0$ in which $F_2(\nu_0(\ell)) =
\varepsilon$, see $(*)_2(f)$.)

In $(*)_3$ we choose
$\cU^{\up}_4,\cU^{\up}_3,g^3_{3,\varepsilon},\alpha^*_3$ and $\langle
s_\delta:\delta \in \cU^{\up}_\ell \rangle$ giving another part of
$\nu_0$.  Then in $(*)_4$ we deal further with $\nu_0$, in particular $s_\delta
\subseteq \kappa_1$ is a stationary subset of
$S^{\kappa_1}_{\kappa_0,j_*}$, promising $F_1(\nu_2(\ell_4)) \in
s_{\delta_2}$.

Next we work on restricting the choices from below, choosing
$\cU^{\dn}_1,\varepsilon^{\dn}_{1,0},\varepsilon^{\dn}_{1,1}$ in
$(*)_5$ promising $\delta_0 \in \cU^{\dn}_1$ so this restricts
$\eta_1$.

Lastly, in $(*)_6$ we choose
$\cU^{\dn}_2,\varepsilon^{\dn}_{2,0},\varepsilon^{\dn}_{2,1}$
promising $\delta_1 \in \cU^{\dn}_2$ (recalling 
$\nu_1 = \rho_{\bar h}(\delta_1,\delta_2)$). 

\begin{claim}  
\label{e4}
Assume $\kappa_1,\kappa_0$ are cardinals and $S$ is a set. 
 There is a function $\bold d:{}^{\omega >}S \rightarrow \bbN$ 
such that $(A) \Rightarrow (B)$ where
\mn
\begin{enumerate}
\item[$(A)$]
\begin{enumerate}
\item[$(a)$]  $F_\iota:S \rightarrow \kappa_\iota$ for $\iota = 0,1$
\sn
\item[$(b)$]   for $\varrho \in {}^{\omega >}S$ and $\iota
  < 2$ we let $F_\iota(\varrho) = \langle F_\iota(\varrho(\ell)):\ell
  < \ell g(\varrho)\rangle$
\sn
\item[$(c)$]  we stipulate $\max \Rang(F_\iota(\langle \rangle)) = -1$
\end{enumerate}
\sn
\item[$(B)$]  $\bold d(\varrho) = \ell^\bullet_4$ \when \,
$\varrho = \eta_0 \char 94 \nu_0 \char 94 \rho \char 94 \nu_1 \char 94
  \eta_1$ satisfies (note that $\ell_1,\ell^\bullet_4 - \ell
  g(\eta_0)$ are places in
$\nu_0,\ell_3$ is a place in $\nu_1,\ell^*_2$ is a place in $\rho$ 
and $\ell^\bullet_2,\ell^\bullet_4$ is a place in $\varrho$ and
  $u \subseteq \{\ell g(\nu_0) + \ell:\ell < \ell g(\nu_0)\}$) the following:
\begin{enumerate}
\item[(a)]
\begin{enumerate}
\item[$(\alpha)$]  $\max \Rang(F_1(\varrho)) =
  \max(\Rang(F_1(\nu_0)) > \max(\Rang(F_1(\eta_0 \char 94 \rho \char 94 \nu_1
  \char 94 \eta_1))$ 
\sn
\item[$(\beta)$]   let $\ell_1 = \min\{\ell <
  \ell g(\nu_0):F_1(\nu_0(\ell)) = \max \Rang(F_1(\varrho))\}$
 so $\ell_1 < \ell g(\nu_0)$
\end{enumerate}
\sn
\item[(b)]
\begin{enumerate}
\item[$(\alpha)$]  $\max \Rang(F_0(\varrho \rest (\ell
  g(\eta_0)+ \ell_1,\ell g(\varrho)))) = \max \Rang(F_0(\rho)) > \\ 
\max \Rang(F_0(\nu_0 \rest [\ell_1,\ell g(\nu_0)) 
\char 94 \nu_1 \char 94 \eta_1)$
\sn
\item[$(\beta)$]  let $\ell^\bullet_2 = \min\{\ell < \ell g(\varrho):\ell \ge
  \ell g(\eta_0) + \ell_1$ and 
$F_0(\varrho(\ell)) = \max \Rang(F_0(\varrho \rest (\ell
  g(\eta_0)+\ell_1,\ell g(\varrho))))\}$ so $\ell_2^\bullet < \ell g(\varrho)$
  and $\ell^*_2 = \ell^\bullet_2 - \ell g(\eta_0 \char 94 \nu_0)$
\sn
\item[$(\gamma)$]  hence $\ell^\bullet_2 \in 
[\ell g(\eta_0 \char 94 \nu_0),\ell g(\eta_0 \char 94 \nu_0 \char 94
\rho))$ and $\ell^*_2 < \ell g(\rho)$
\end{enumerate}
\sn
\item[(c)]
\begin{enumerate}
\item[$(\alpha)$]  $\max \Rang(F_1(\nu_0)) > \max
  \Rang(F_1(\varrho \rest [\ell^\bullet_2,\ell g(\varrho))) 
= \max \Rang(F_1(\nu_1)) > \max\{F_1(\rho(\ell)):\ell \in
[\ell^*_2,\ell g(\rho))\}$
\sn
\item[$(\beta)$]  $\ell_3$ is such that
\sn
\begin{enumerate}
\item[$\bullet_1$]  $\ell_3 < \ell g(\nu_1)$
\sn
\item[$\bullet_2$]  $F_1(\nu_1(\ell_3)) = 
\max\{F_1(\varrho)(\ell):\ell \ge \ell^\bullet_2\}$
\sn
\item[$\bullet_3$]   $\ell_3$ is minimal under the above
\end{enumerate}
\end{enumerate}
\sn
\item[(d)]
\begin{enumerate}
\item[$(\alpha)$]   let $u := \{\ell:\ell \le \ell^\bullet_2$ and 
$F_1(\varrho)(\ell) \ge F_1(\nu_1(\ell_3))\}$ 
\sn
\item[$(\beta)$]  $\ell^\bullet_4 \in u$ is such that
\sn
\begin{enumerate}
\item[$\bullet_1$]  $F_1(\varrho(\ell^\bullet_4)) =
  \min\{F_1(\varrho(\ell)):\ell \in u\}$
\sn
\item[$\bullet_2$]   under $\bullet_1,\ell^\bullet_4$ is minimal
\sn
\item[$\bullet_3$]   \underline{notation}: if $\ell^\bullet_4 \in 
[\ell g(\eta_0),\ell g(\eta_0 \char 94 \nu_0))$ then we let $\ell^*_4 =
\ell^\bullet_4 - \ell g(\eta_0)$.
\end{enumerate}
\end{enumerate}
\end{enumerate}
\end{enumerate}
\end{claim}

\begin{PROOF}{\ref{e4}}
Assume $\varrho \in {}^{\omega >} S$.  We have to show that $\bold d$
is well defined, i.e. $\bold d(\varrho) = \ell^\bullet_4$
does not depend on the specific representation of $\varrho$ as $\eta_0
\char 94 \nu_0 \char 94 \rho \char 94 \nu_1 \char 94 \eta_1$, i.e.
we shall prove that $\ell^\bullet_4$ depends on $\varrho$ only.

Toward this
\mn
\begin{enumerate}
\item[$(a)$]  $\ell g(\eta_0) + \ell_1$ depends on $\varrho$ only
\end{enumerate}
\mn
[Why?  Let $\ell^\bullet_1$ be the first natural number so that
$F_1(\varrho(\ell^\bullet_1)) = \max\Rang(F_1(\varrho))$.  By the strict $>$
in $(B)(a)(\alpha)$ we must have $\ell g(\eta_0) \le \ell^\bullet_1$.  Although
one can decompose $\varrho$ in different ways, yielding different
values to $\ell g(\eta_0)$, the sum $\ell g(\eta_0) + \ell_1$ will be
always $\ell^\bullet_1$, by the definition of $\ell_1$.  Now since only
$\varrho$ is mentioned in the definition of $\ell^\bullet_1$ we conclude
that $\ell g(\eta_0) +\ell_1 = \ell^\bullet_1$ depends on $\varrho$ only.]
\mn
\begin{enumerate}
\item[$(b)$]  $\ell^\bullet_2$ depends 
on $\varrho$ only by a similar argument, this time for the function $F_0$
\sn
\item[$(c)$]  $\ell g(\eta_0 \char 94 \nu_0 \char 94 \rho) + \ell_3$ 
depends on $\varrho$ only (for this statement notice that $\rho \ne
\langle \rangle$, by $(b)(\alpha)$)
\sn
\item[$(d)$]  $\{\ell g(\eta_0) + \ell:\ell \in u\}$ depends on $\varrho$ only
\sn
\item[$(e)$]  $\ell^\bullet_4$ depends on $\varrho$ only.
\end{enumerate}
\mn
By (e) clearly we are done.
\end{PROOF}

\begin{theorem}  
\label{e8}
Assume $\aleph_0 \le \theta = \cf(\theta),\lambda \ge \theta^{+3}$ and
$\lambda$ is a successor of a regular cardinal.  
\Then \, $\Pr_1(\lambda,\lambda,\lambda,\theta)$ holds.
\end{theorem}

\begin{PROOF}{\ref{e8}}
Firstly, let us spell out the definition of $\Pr_1$.

Recall that $\lambda \ge \mu \ge \sigma,\theta_0,\theta_1$ and let
$\bar\theta = (\theta_0,\theta_1)$. 
$\Pr_1(\lambda,\mu,\sigma,\bar\theta)$ means
that there exists a function $\bold c:[\lambda]^2 \rightarrow \sigma$
such that for every two disjoint sequences $\langle
\zeta^0_{\alpha,i}:\alpha < \mu,i < \bold i_0\rangle,\langle
\zeta^1_{\alpha,i}:\alpha < \mu,i < \bold i_1\rangle$ of ordinals $<
\lambda$ (without repetitions) such that $\bold i_0 < \theta_0,\bold
i_1 < \theta_1$ and for every $\gamma < \sigma$, one can find
$\alpha_0 < \alpha_1 < \mu$ so that:
\mn
\begin{enumerate}
\item[$(*)$]  if $i_0 < \bold i_0$ and $i_1 < \bold i_1$ then $\bold
  c(\zeta^0_{\alpha_0,i_0},\zeta^1_{\alpha_1,i_1}) = \gamma$.
\end{enumerate}
\mn
It follows from the definition that if $\theta'_1 \le \theta_1$ and
$\Pr_1(\lambda,\mu,\sigma,(\theta_0,\theta_1))$ then $\Pr_1
(\lambda,\mu,\sigma,(\theta_0,\theta'_1))$.  Let $\theta_0 =
\theta,\theta_1 =\theta^+$ by Theorem \ref{e16}
below we have $\Pr_1(\lambda,\lambda,\lambda,(\theta_0,\theta_1))$
and since $\theta_0 < \theta_1$ we have by the previous sentence
$\Pr_1(\lambda,\lambda,\lambda,(\theta_0,\theta_0))$ which is also denoted
$\Pr_1(\lambda,\lambda,\lambda,\theta_0)$, see
Observation \ref{d36}, so we are done by noticing that $\theta_0$ of
\ref{e16} is $\theta$ here.
\end{PROOF}

\begin{remark}  
\label{e9}
1) Can we replace $\theta$ by $(\theta^+,\theta)$? 

\noindent
2) Or, at least when $\theta = \aleph_0,\lambda = \aleph_2$ 
for $(\theta,\theta^+)$ with an ultrafilter on the 
$< \theta^+$ sets? and 2 colours? may try to use the proof of the
$\aleph_2$-c.c. not productive from \cite{Sh:572}.

\noindent
3) For many purposes, $\Pr_1(\lambda,\lambda,2,(\theta,\theta^+))$
suffices and for this the proof (in \ref{e16}) is somewhat simpler.
\end{remark}

\begin{conclusion}
\label{e12}
Assume $\lambda = \partial^+,\partial = \cf(\partial) > \theta^+,
\theta = \cf(\theta) \ge \aleph_0$
\mn
\begin{enumerate}
\item[(a)]  if there is $\chi = \chi^{< \theta} < \lambda \le 2^\chi$ 
and $\chi \ge \sigma$ (so $\sigma \le \partial$), \then \,
$\Pr_0(\lambda,\lambda,\sigma,\theta)$
\sn
\item[(b)]  if $\chi = \partial$ satisifes $\chi = \chi^{< \theta}$
  then $\Pr_0(\lambda,\lambda,\lambda,\theta)$.
\end{enumerate}
\end{conclusion}

\begin{PROOF}{\ref{e12}}
\underline{Clause (a)}:

We apply \ref{d43}(4) with
$(\lambda,\lambda,\chi,\sigma,\theta,\theta)$ here standing for
$(\lambda,\mu,\chi,\sigma,\theta_0,\theta_1)$
there.  We have to check the assumption of \ref{d43}(4), the main
point is ``$\Pr_1(\lambda,\lambda,\sigma,(\theta,\theta))"$ which holds by
Theorem \ref{e8}, the other assumptions are straightforward hence we
get the conclusion, i.e. $\Pr_0(\lambda,\lambda,\sigma,\bar\theta)$.
\medskip

\noindent
\underline{Clause (b)}:

First, $\Pr_0(\lambda,\lambda,\partial,\theta)$ holds as we can apply
Clause (a) with $(\lambda,\partial,\partial,\partial,\theta)$ here
standing for $(\lambda,\partial,\chi,\sigma,\bar\theta)$ there.

Second, we get $\Pr_0(\lambda,\lambda,\lambda,\theta)$ holds as we can apply
\ref{d43}(3) with $(\lambda,\partial,\theta)$ here standing for
$(\lambda,\sigma,\theta)$ there.
\end{PROOF}

\begin{theorem}
\label{e16}
If $\lambda$ is a successor of a regular cardinal, $\lambda \ge
\theta^+_1$ and $\theta_1 > \theta_0 \ge \aleph_0$ are 
regular cardinals, \then \,
$\Pr_{1}(\lambda,\lambda,\lambda,(\theta_0,\theta_1))$.
\end{theorem}

\begin{PROOF}{\ref{e8}}
\medskip

\noindent
\underline{Stage A}:

Let $\partial$ be the regular cardinal such that $\lambda
= \partial^+$, so $\partial > \theta_1$.

Below we shall choose $\sigma$ and $\kappa_\iota$ (for $\iota=0,1,2$)
to help in using this proof for proving other theorems.

Let $\sigma = \lambda$.
Let $S \subseteq S^\lambda_\partial$ be stationary and $h:\lambda
\rightarrow \sigma$ be such that $\alpha < \lambda \Rightarrow
h(\alpha) < 1 + \alpha,h \rest (\lambda \backslash S)$ is
constantly zero and $S^*_\gamma := \{\delta \in S:h(\delta) = \gamma\}$
is a stationary subset of $\lambda$ for every $\gamma < \lambda$.  Let
$(\kappa_0,\kappa_1,\kappa_2) = (\theta_0,\theta_1,\sigma)$ and let
$F_\iota:\lambda = \sigma \rightarrow \kappa_\iota$ for $\iota = 0,1,2$ be such
that for every $(\varepsilon_0,\varepsilon_1,\varepsilon_2) \in
(\kappa_0 \times \kappa_1 \times \kappa_2)$ the set
$W_{\varepsilon_0,\varepsilon_1,\varepsilon_2}(\kappa) = \{\gamma \in
S^\lambda_\kappa:F_\iota(\gamma) = \varepsilon_\iota$ for $\iota \le 2\}$ is a
stationary subset of $\lambda$ for every $\kappa = \cf(\kappa) < \lambda$.

Let $\bar e = \langle e_\alpha:\alpha < \lambda\rangle$ be such that
\mn
\begin{enumerate}
\item[$\odot_1$]
\begin{enumerate}
\item[(a)]   if $\alpha = 0$ then $e_\alpha = \emptyset$
\sn
\item[(b)]   if $\alpha = \beta +1$ then $e_\alpha = \{\beta\}$
\sn
\item[(c)]   if $\alpha$ is a limit ordinal then $e_\alpha$ is a club of
  $\alpha$ of order type $\cf(\alpha)$ disjoint
to $S^\lambda_\partial$ hence to $S$.
\end{enumerate}
\end{enumerate}
\mn
Let\footnote{For successor of regular we can omit $h_\alpha$ and below
  replace $\bar h$ and $h^-$ by $h$ and even let $\rho_h(\beta,\alpha)
  = \langle h(\gamma_\ell(\beta,\alpha)):\ell <
  k(\beta,\alpha)\rangle$; but for other cases the present version is
  better, see more \cite[Ch.III,\S4]{Sh:g}.  But in later stages we may
  use $h$ directly, e.g. the proof of $(*)_1$.}
$h_\alpha = h \rest e_\alpha$ for $\alpha < \lambda$ and $\bar h =
\langle h_\alpha:\alpha < \lambda\rangle$. Note that $h_\alpha$ is
non-zero only for successor $\alpha$.  We shall mostly use the
$h_\alpha$'s rather than $h$.

Now (using $\bar e$) for $0 < \alpha < \beta < \lambda$, let

\[
\gamma(\beta,\alpha) := \min\{\gamma \in e_{\beta}:\gamma \ge
\alpha\}.
\]

\mn
Let us define $\gamma_{\ell}(\beta,\alpha)$:

\[
\gamma_{0}(\beta,\alpha) = \beta,
\]

\[
\gamma_{\ell +1}(\beta,\alpha) = \gamma(\gamma_{\ell}(\beta,\alpha),
\alpha) \text{ (if defined)}.
\]

\mn
If $0 < \alpha < \beta < \lambda$, let $k(\beta,\alpha)$ be the maximal $k < 
\omega$ such that $\gamma_{k}(\beta,\alpha)$ is defined (equivalently
is equal to $\alpha$) and let
$\rho_{\beta,\alpha} = \rho(\beta,\alpha)$  be the sequence  

\[
\langle \gamma _{0}(\beta,\alpha),\gamma_{1}(\beta,\alpha),\ldots,
\gamma_{k(\beta,\alpha)-1}(\beta,\alpha)\rangle.
\]

\mn
Let $\gamma_{\ell t}(\beta,\alpha) =
\gamma_{k(\beta,\alpha)-1}(\beta,\alpha)$ where $\ell t$ stands for last.   

Let

\[
\rho_{\bar h}(\beta,\alpha) = \langle h_{\gamma_{\ell}(\beta,\alpha)}
(\gamma_{\ell +1}(\beta,\alpha)): \ell < k(\beta,\alpha)\rangle
\]

\mn
and we let $\rho(\alpha,\alpha)$ and  $\rho_{\bar h}(\alpha,\alpha)$ 
be the empty sequence.

\noindent
Now clearly:
\mn 
\begin{enumerate}
\item[$\odot_2$]   if $0 < \alpha < \beta < \lambda$ then $\alpha \le
\gamma(\beta,\alpha) < \beta$
\end{enumerate}
\mn
hence
\mn 
\begin{enumerate}
\item[$\odot_3$]   if $0 < \alpha < \beta < \lambda,0 < \ell < \omega$, and
$\gamma_{\ell}(\beta,\alpha)$ is well defined, then

\[
\alpha \le \gamma_{\ell}(\beta,\alpha) < \beta 
\]
\end{enumerate}
\mn
and
\mn 
\begin{enumerate}
\item[$\odot_4$]  if $0 < \alpha < \beta < \lambda$, then $k(\beta,\alpha)$ is
well defined and letting $\gamma_{\ell} :=
\gamma_{\ell}(\beta,\alpha)$ for $\ell \le k(\beta,\alpha)$ we have
 
\[
\alpha = \gamma_{k(\beta,\alpha)} < \gamma_{\ell t}(\beta,\alpha) = 
\gamma_{k(\beta,\alpha)-1} < \cdot \cdot \cdot < \gamma_{1} < \gamma_{0} = 
\beta 
\]

\[
\text{ and } \alpha \in e_{\gamma_{\lt}(\beta,\alpha)}
\]

\mn
i.e. $\rho(\beta,\alpha)$ is a (strictly) decreasing finite 
sequence of ordinals, starting with $\beta$, ending with 
$\gamma_{\ell t}(\beta,\alpha)$ of length $k(\beta,\alpha)$.
\end{enumerate}
\mn
Note that if $\alpha \in S,\alpha < \beta$ then
$\gamma_{\lt}(\beta,\alpha) = \alpha +1$.

Also
\mn 
\begin{enumerate}
\item[$\odot_5$]  if $\delta$ is a limit ordinal and 
$\delta < \beta < \lambda$,
\then \, for some $\alpha_{0} < \delta$ we have: $\alpha_{0} \le \alpha <
\delta$ \underline{implies}:
\sn
\begin{enumerate}
\item[$(i)$]  for $\ell < k(\beta,\delta)$ we have 
$\gamma_{\ell}(\beta,\delta) = \gamma_{\ell}(\beta,\alpha)$   
\sn
\item[$(ii)$]  $\delta \in \nacc(e_{\gamma_{\ell t}(\beta,\delta)})
\Leftrightarrow \delta = \gamma_{k(\beta,\delta)}(\beta,\delta) =
\gamma_{k(\beta,\delta)}(\beta,\alpha) \Leftrightarrow 
\neg[\gamma_{k(\beta,\delta)}(\beta,\delta) = \delta > 
\gamma_{k(\beta,\delta)}(\beta,\alpha)]$   
\sn
\item[$(iii)$]  $\rho(\beta,\delta) \trianglelefteq
  \rho(\beta,\alpha)$; i.e. is an initial segment
\sn
\item[$(iv)$]  $\delta \in \nacc(e_{\gamma_{\ell t}(\beta,\delta)})$
  (here always holds if $\delta \in S$) implies:
\begin{itemize}
\item  $\rho(\beta,\delta) \char 94 \langle \delta \rangle \trianglelefteq 
\rho(\beta,\alpha)$ hence 
\sn
\item  $\rho_{\bar h}(\beta,\delta) \char 94 
\langle h_{\gamma_{\ell t}(\beta,\delta)}
(\delta )\rangle \trianglelefteq \rho_{\bar h}(\beta,\alpha)$.
\end{itemize}
\sn
\item[$(v)$]  if $\cf(\delta) = \partial$ \then \, we have 
$\gamma_{\ell t}(\beta,\delta) = \delta +1$
\sn
\item[$(vi)$]  if $\cf(\delta) = \partial$ and $\delta \in e_\alpha$, 
\then \, necessarily $\alpha = \delta +1$.
\end{enumerate}
\end{enumerate}
\mn
Why?  Just let

\[
\alpha_{0} = \Max\{\sup(e_{\gamma_{\ell}(\beta,\delta)} \cap \delta) +
1: \ell < k(\beta,\delta) \text{ and } \delta \notin
\acc(e_{\gamma_{\ell}(\beta,\delta)})\}.
\]

\mn
Notice that if $\ell < k(\beta,\delta)-1$ then $\delta \notin
\acc(e_{\gamma_\ell(\beta,\delta)})$ is immediate.

Note that the outer maximum (in the choice of $\alpha_0$) 
is well defined as it is over a finite non-empty set of 
ordinals.  The inner $\sup$ is on the empty set (in which case we get
zero) or is the maximum
 (which is well defined) as $e_{\gamma_{\ell}(\beta,\delta)}$ is a 
closed subset of $\gamma_{\ell}(\beta,\delta),\delta < \gamma_{\ell}
(\beta,\delta)$ and $\delta \notin
\acc(e_{\gamma_{\ell}(\beta,\delta)})$ - as this is required.  For
clauses (v), (vi) recall $\delta \in S^\lambda_\partial$ and $e_\gamma \cap
S^\lambda_\partial = \emptyset$ when $\gamma$ is a limit ordinal and
$e_\gamma = \{\gamma-1\}$ when $\gamma$ is a successor ordinal.
\mn
\begin{enumerate}
\item[$\odot_6$]  
\begin{enumerate}
\item[(a)]  if $0 < \alpha < \beta < \lambda,
\ell < k(\beta,\alpha),\gamma =
\gamma_{\ell}(\beta,\alpha)$ \then \, $\rho(\beta,\alpha) =
\rho(\beta,\gamma) \char 94 \rho(\gamma,\alpha)$ and 
$\rho_{\bar h}(\beta,\alpha) = 
\rho_{\bar h}(\beta,\gamma) \char 94 \rho_{\bar h}(\gamma,\alpha)$
\sn
\item[(b)]  if $0 < \alpha_0 < \ldots < \alpha_k$ and
  $\rho(\alpha_k,\alpha_0) = \rho(\alpha_k,\alpha_{k-1}) \char 94
  \ldots \char 94 \rho(\alpha_1,\alpha_0)$ \then \, this holds for any
  subsequence of $\langle \alpha_0,\dotsc,\alpha_k\rangle$.
\end{enumerate}
\end{enumerate}
\mn
Now apply Claim \ref{e4} with $\lambda,\kappa_1,\kappa_0,F_1,F_0$ here
standing for $S,\kappa_1,\kappa_0,F_1,F_0$ there and get 
$\bold d:{}^{\omega >} \lambda \rightarrow \bbN$.

Lastly, we define the colouring; as the proof is somewhat simpler if we
use only $\kappa_1$ colours (which suffice for many
purposes) we define two colourings: $\bold c_1$ with $\kappa_1$ 
colours and $\bold c_2$ with $\kappa_2 = \lambda$ colours,
as follows:
\mn
\begin{enumerate}
\item[$\odot_7$]
\begin{enumerate}
\item[(a)]   choose a function $h':\kappa_1 \rightarrow
  \kappa_1$ such that $S^{\kappa_1}_{\kappa_0,\varepsilon} :=
  \{\delta \in S^{\kappa_1}_{\kappa_0}:h'(\delta) = \varepsilon\}$
is stationary in $\kappa_1$ for every $\varepsilon < \kappa_1$
\sn
\item[(b)]  if $\eta = \langle \zeta_0,\dotsc,\zeta_{n-1}\rangle$ then
we let $h'(\eta) = \langle h'(\zeta_0),\dotsc,h'(\zeta_{n-1})\rangle$
\sn
\item[(c)]  $\bold c_1:[\lambda]^2 \rightarrow \kappa_1$ is defined 
for $\alpha < \beta$ by $\bold c_1(\{\alpha,\beta\}) = h'(F_1(\rho_{\bar h}
(\beta,\alpha)))(\ell^1_{\beta,\alpha})$ 
 where
$\ell^1_{\beta,\alpha} = \bold d(\rho_{\bar h}(\beta,\alpha))$.
\end{enumerate}
\end{enumerate}
\mn
Clearly
\mn
\begin{enumerate}
\item[$\odot_8$]  we can demand on $h'_1$ that we can choose $h'_2$ such that:
\sn
\begin{enumerate}
\item[(a)]  $h'_1,h'_2$ are functions with domain $\kappa_1$
\sn
\item[(b)]  $h'_1$ is onto $\kappa_1$
\sn
\item[(c)]  $h'_2$ is onto $\bbN$
\sn
\item[(d)]  for every $\zeta < \kappa_1$ and $n <
  \omega$ the set $S_{\kappa_1,\zeta,n} = \{\varepsilon < \kappa_1:
h'_1(\varepsilon) = \zeta$ and $h'_2(\varepsilon) =n\}$ is stationary
\end{enumerate}
\sn
\item[$\odot_9$]  the colouring $\bold c_2$ with $\lambda$ colours is
  chosen as follows: for 
$\alpha < \beta < \lambda,\bold c_2(\{\alpha,\beta\}) =
(F_2(\rho_{\bar h}(\beta,\alpha)))(\ell^2_{\beta,\alpha})$ where
letting $\varepsilon_{\alpha,\beta} = \bold c_1(\{\alpha,\beta\})$ we have
$\ell^2_{\beta,\alpha}$ is the $h'_2(\varepsilon_{\beta,\alpha})$-th
member of the\footnote{So $\bold d$ is used only via the definition of
  $\ell^2_{\beta,\alpha}$.}  set $\{\ell < \ell g(\rho_{\bar
  h}(\beta,\alpha)):F_1(\rho_{\bar h}(\beta,\alpha))(\ell) = 
h'_1(\varepsilon_{\beta,\alpha})\}$ if this set has 
$> h'_2(\varp_{\alpha,\beta})$ members and is zero otherwise. 
\end{enumerate}
\mn
\underline{Stage B}:

So we have to prove that the colouring $\bold c = \bold c_1$ (with
$\kappa_1$ colours) and moreover $\bold c = \bold c_2$ (with $\lambda$
colours) is as required.  

Now for the rest of the proof assume:
\mn
\begin{enumerate}
\item[$\boxplus$]
\begin{enumerate}
\item[(a)]  $t_{\alpha} \subseteq \lambda$ for every $\alpha < \lambda$
\sn
\item[(b)]  $t_\alpha = t^0_\alpha \cup t^1_\alpha$ and 
$1 \le |t^\iota_\alpha| < \theta_\iota$ for $\iota < 2$
\sn
\item[(c)]  $\alpha \ne \beta \Rightarrow 
t_{\alpha} \cap t_{\beta}  = \emptyset$
\sn
\item[(d)]  $j_* < \kappa_1$ (when dealing with $\bold c_1$) 
or $j_* < \sigma$ (when dealing with $\bold c_2$).
\end{enumerate}
\end{enumerate}
\mn
Clearly (by $\boxplus(c)$), we can choose $\beta_\alpha$ by induction on
$\alpha < \lambda$ by $\beta_\alpha = \min\{\beta:\beta > \alpha$ and
$\min(t_\beta) > \alpha + \sup(\cup\{t_{\beta_{\alpha(1)}}:\alpha(1) <
\alpha\})\}$.  Now can use $t'_\alpha = t_{\beta_\alpha}$ for $\alpha
< \lambda$, hence:
\mn
\begin{enumerate}
\item[$(*)_0$]   \wilog \, $\alpha < \min(t_\alpha)$ and $\alpha <
  \beta \Rightarrow \sup(t_\alpha) < \min(t_\beta)$.
\end{enumerate}
\mn
We have to prove that for some $\alpha_0 < \alpha_1 < \lambda$ for
every $(\zeta_0,\zeta_1) \in t^0_{\alpha_0} \times t^1_{\alpha_1}$ we have
$\bold c\{\zeta_0,\zeta_1\} = j_*$.
\mn
\begin{enumerate}
\item[$(*)_1$]  We can find $\cU^{\up}_1,\alpha^*_1,
\varepsilon^{\up}_{1,1}$ such that:
\sn
\begin{enumerate}
\item[$(a)$]  $\cU^{\up}_1 \subseteq S$ is stationary
\sn
\item[$(b)$]  $h \rest \cU^{\up}_1$ is constantly 0 (so actually
  $\cU^{\up}_1 \subseteq S^*_0$)
\sn
\item[$(c)$]  $\alpha^*_1 < \min(\cU^{\up}_1)$ and
$\varepsilon^{\up}_{1,1} < \kappa_1$
\sn
\item[$(d)$]  if $\delta \in \cU^{\up}_1$ and $\alpha \in
[\alpha^*_1,\delta),\beta \in t^1_\delta$ (treating $t^0_\delta$
is unreasonable because $t^1_\delta$ may be of cardinality $\ge
  \theta_0 = \kappa_1,\varp_{1,0}$ is defined for notational simplicity)
\then \,:
\sn
\begin{itemize}
\item  $\rho_{\beta,\delta} \char 94 \langle \delta \rangle
  \trianglelefteq \rho_{\beta,\alpha}$
\sn
\item  $\Rang(F_1(\rho_{\bar h}(\beta,\delta))) 
\subseteq \varepsilon^{\up}_{1,1}$.
\end{itemize}
\end{enumerate}
\end{enumerate}
\mn
[Why?  For every $\delta \in S^*_0 \subseteq S$ and $\zeta \in t_\delta$ let
$\alpha^*_{1,\delta,\zeta} < \delta$ be such that $(\forall \alpha)
(\alpha \in [\alpha^*_{1,\delta,\zeta},\delta) \Rightarrow
\rho_{\zeta,\delta} \char 94 \langle \delta \rangle 
\trianglelefteq \rho_{\zeta,\alpha})$, it exists by $\odot_5$ of Stage A.

Let $\alpha^*_{1,\delta} = \sup\{\alpha^*_{1,\delta,\zeta}:\zeta \in
t_\delta\}$ and for $\iota=1$ let $\varepsilon^{\up}_{1,1,\delta} =
\sup\{F_1(h(\gamma_\ell(\zeta,\delta)))+1:\zeta \in t^1_\delta$ and $\ell
< k(\zeta,\delta)\} = \sup \cup\{\Rang(F_1(\rho_{\bar h}
(\beta,\delta)) +1):\beta \in t^1_\delta\})$; 
as $\cf(\delta) = \partial = \cf(\partial) > |t^1_\delta|$
and $\kappa_1 = \cf(\kappa_1) \ge \theta_1 > |t^1_\delta|$, necessarily
$\alpha^*_{1,\delta} < \delta$ and $\varepsilon^{\up}_{1,1,\delta}
< \kappa_\iota$.

Lastly, there are $\alpha^*_1 < \lambda$ and
$\varepsilon^{\up}_{1,0} < \kappa_0,\varepsilon^{\up}_{1,1} <
\kappa_1$ and $\cU^{\up}_1 \subseteq S^*_0$ as required in $(*)_1$ by 
Fodor lemma.]
\mn
\begin{enumerate}
\item[$(*)_2$]  for each $\varepsilon \in S^{\kappa_1}_{\kappa_0}
  \backslash \varepsilon^{\up}_{1,1}$ we can find
$g_{2,\varepsilon},\cU^{\up}_{2,\varepsilon},
\gamma^*_\varepsilon,\alpha^*_{2,\varepsilon},\ell_{2,\varepsilon}$ such that:
\sn
\begin{enumerate}
\item[$(a)$]  $\gamma^*_\varepsilon < \lambda$ satisfies
  $F_2(\gamma^*_\varepsilon) = j_*,F_1(\gamma^*_\varepsilon) =
  \varepsilon,F_0(\gamma^*_\varepsilon) = 0$
\sn
\item[$(b)$]  $\cU^{\up}_{2,\varepsilon} \subseteq
  S^*_{\gamma^*_\varepsilon}$ is stationary
\sn
\item[$(c)$]  $\alpha^*_1 < \alpha^*_{2,\varepsilon} <
  \min(\cU^{\up}_{2,\varepsilon})$ 
\sn
\item[$(d)$]   $g_{2,\varepsilon}$ is a function with domain 
$\cU^{\up}_{2,\varepsilon}$ such that $\delta \in \cU^{\up}_{2,\varepsilon}
  \Rightarrow \delta < g_{2,\varepsilon}(\delta) \in \cU^{\up}_1$
\sn
\item[$(e)$]  if $\delta \in \cU^{\up}_{2,\varepsilon}$ and $\alpha \in
[\alpha^*_{2,\varepsilon},\delta)$ and $\beta \in
  t_{g_{2,\varepsilon}(\delta)}$ then
$\rho_{g_{2,\varepsilon}(\delta),\delta} \char 94 \langle \delta \rangle
  \trianglelefteq \rho_{g_{2,\varepsilon}(\delta),\alpha}$ hence
  (recalling $\odot_6,(*)_1(d)$)
\sn
\begin{itemize}
\item  if $\beta \in t_{g_{2,\varp}(\delta)}$ then 
$\rho_{\beta,\delta} \char 94 
\langle \delta \rangle \trianglelefteq \rho_{\beta,\alpha}$
\end{itemize}
\sn
\item[$(f)$]  $\ell^*_{2,\varepsilon}$ is well defined where for any
  $\delta \in \cU^{\up}_{2,\varepsilon}$ we have
\newline
$\ell^*_{2,\varepsilon} = \ell
g(\rho_{g_{2,\varepsilon}(\delta),\delta})$ hence if $\alpha \in
(\alpha^*_{2,\varp},\delta)$ then
$\rho_{g_{2,\varp}(\delta),\alpha}(\ell^*_{2,\varp}) = \delta$.  
\sn
\item[$(g)$]  Lastly, if $\alpha \in
(\alpha^*_{2,\varp},\delta)$ then $\ell^\bullet_{2,\varepsilon} 
= \min\{\ell:\ell < 
\ell g(\rho_{g_{2,\varepsilon}(\delta),\alpha})$ and
$F_1(\rho_{\bar h}(g_{2,\varepsilon}(\delta),\alpha))(\ell) =
\varepsilon\}$ so $\ell^\bullet_{2,\varp} \le \ell^*_{2,\varp}$; 
recall that $\varepsilon > \varepsilon^{\up}_{1,1}$ hence
 necessarily $\beta \in t_{g_{2,\varepsilon}(\delta)} \Rightarrow
\varepsilon > \sup\Rang(F_1(\rho_{\bar h}(\beta,g_{2,\varepsilon}(\delta))))$. 
\end{enumerate}
\end{enumerate}
\mn
[Why?  First, choose $\gamma^*_\varepsilon$ as in clause (a) of $(*)_2$,
(possible by the choice of $F_0,F_1,F_2$ in the beginning of Stage A).  
Second, define $g'_\varepsilon:S^*_{\gamma^*_\varepsilon} 
\rightarrow \cU^{\up}_1$ such that $\delta \in
S^*_{\gamma^*_\varepsilon} \Rightarrow \delta <
  g'_\varepsilon(\delta) \in \cU^{\up}_1$.  
Third, do as in the proof of $(*)_1$ 
above for each $\delta \in S^*_{\gamma^*_\varepsilon}$
separately, i.e. find $\alpha'_{2,\varepsilon,\delta} < \delta$ 
above $\alpha^*_1$ and
$\ell^*_{2,\varepsilon,\delta},\ell^\bullet_{2,\varp,\delta}$ 
such that the parallel of clauses $(c),(e),(f),(g)$ of 
$(*)_2$ holds.  Fourth, use Fodor lemma to get a
stationary $\cU^{\up}_{2,\varepsilon} \subseteq
 S^*_{\gamma^*_\varepsilon}$ such that $\langle
(\alpha'_{2,\varepsilon,\delta},\ell^*_{2,\varepsilon,\delta},\ell^\bullet_{2,\varp,\delta}):\delta \in \cU^{\up}_{2,\varepsilon}\rangle$ is constantly 
$(\alpha^*_{2,\varepsilon},\ell^*_{2,\varepsilon},\ell^\bullet_{2,\varp})$
and lastly let
$g_{2,\varepsilon} = g'_\varepsilon \rest \cU^{\up}_{2,\varepsilon}$.]
\mn
\begin{enumerate}
\item[$(*)_3$]  we can find $\cU^{\up}_3,\bar g^3,\alpha^*_3$ such that:
\sn
\begin{enumerate}
\item[$(a)$]  $\cU^{\up}_3 \subseteq S$ is stationary 
\sn
\item[$(b)$]  $\min(\cU^{\up}_3) > \alpha^*_3 >
 \sup\{\alpha^*_{2,\varepsilon}:\varepsilon \in
 S^{\kappa_1}_{\kappa_0} \backslash \varp^{\up}_{1,1}\}$
\sn
\item[$(c)$]  $\bar g^3 = \langle g_{3,\varepsilon}:\varepsilon \in
  S^{\kappa_1}_{\kappa_0} \backslash \varepsilon^{\up}_{1,1}\rangle$
\sn
\item[$(d)$]  $g_{3,\varepsilon}$ is a function with domain $\cU^{\up}_3$
\sn
\item[$(e)$]  if $\delta \in \cU^{\up}_3$ and $\varp \in
  S^{\kappa_1}_{\kappa_0} \backslash \varp^{\up}_{1,1}$ then $\delta <
  g_{3,\varepsilon}(\delta) \in \cU^{\up}_{2,\varepsilon}$
\sn
\item[$(f)$]  if $\alpha \in [\alpha^*_3,\delta),\delta \in
\cU^{\up}_3$ and $\varepsilon \in S^{\kappa_1}_{\kappa_0} \backslash
\varepsilon^{\up}_{1,1}$ then $\rho_{g_{3,\varepsilon}(\delta),\delta} 
\char 94 \langle \delta \rangle \trianglelefteq 
\rho_{g_{3,\varepsilon}(\delta),\alpha}$ hence
\sn
\item[$(f)'$]  if in addition $\beta \in
  t^1_{g_{2,\varepsilon}(g_{3,\varepsilon}(\delta))}$ 
then $\rho_{\beta,\delta} \char 94 \langle \delta 
\rangle \trianglelefteq \rho_{\beta,\alpha}$ this follows.
\end{enumerate}
\end{enumerate}
\mn
[Why?  First, let $\alpha^*_2 = \sup\{\alpha^*_{2,\varp} +1:\varp \in
S^{\kappa_1}_{\kappa_0} \backslash \varp^{\up}_{1,1}\} < \lambda$ and
 choose $g''_\varepsilon:S \backslash \alpha^*_2 \rightarrow
  \cU^{\up}_{2,\varepsilon}$ such that $g''_\varepsilon(\delta) >
  \delta$ for every $\delta \in S \backslash \alpha^*_2$ 
and second for each $\delta \in S \backslash \alpha^*_2$ choose 
$\alpha^*_{3,\delta} < \delta$ as in clauses (f),(f)$'$ 
of $(*)_3$, i.e. such that $\alpha \in
[\alpha^*_{3,\delta},\delta) \Rightarrow
 \rho_{g''_\varepsilon(\delta),\delta} \char 94 \langle \delta \rangle
 \trianglelefteq \rho_{g''_\varepsilon(\delta),\alpha}$ for
  every $\varepsilon \in S^{\kappa_1}_{\kappa_0} \backslash
\varepsilon^\up_{1,1}$ and such that the relevant part of clause (b)
  of $(*)_3$, holds, that is, $\alpha^*_{3,\delta} > \alpha^*_2 =
\sup\{\alpha^*_{2,\varepsilon}:\varepsilon < S^{\kappa_1}_{\kappa_0}
\backslash \varepsilon^\up_{1,1}\}$, possible as $\kappa_1 < \partial$.  
Third, use Fodor lemma to find $\alpha^*_3 <
  \lambda$ such that $\cU^{\up}_3 = \{\delta \in S:\alpha^*_{3,\delta}
  = \alpha^*_3\}$ is a stationary subset of $\lambda$.  Fourth, let
  $g_{3,\varepsilon} = g''_\varepsilon \rest \cU^{\up}_3$.] 
\mn
\begin{enumerate}
\item[$(*)_4$]  recalling\footnote{Recall that in this stage we are dealing
  with $\bold c = \bold c_1$ hence $j_* < \kappa_1$.} 
$j_* < \kappa_1$, there are
$\cU^{\up}_4,\varepsilon^*_{4,1},\varepsilon^*_{4,0}$ 
and $\langle s_\delta:\delta \in \cU^{\up}_4\rangle$ such that:
\sn
\begin{enumerate}
\item[$(a)$]  $\cU^{\up}_4 \subseteq \cU^{\up}_3$ is a stationary
  subset of $\lambda$
\sn
\item[$(b)$]  $\varepsilon^{\up}_{1,1} < \varepsilon^{\up}_{4,1} < \kappa_1$
  and $\varepsilon^{\up}_{4,0} < \kappa_0$
\sn
\item[$(c)$]  if $\delta \in \cU^{\up}_4$ then $s_\delta$ is a stationary (in
 $\kappa_1$) subset of $S^{\kappa_1}_{\kappa_0,j_*} 
\backslash \varepsilon^\up_{4,1}$
\sn
\item[$(d)$]  if $\delta \in \cU^{\up}_4,\varepsilon \in s_\delta$
  \then \,
\begin{enumerate}
\item[$(\alpha)$]  $\Rang(F_1(\rho_{\bar h}
(g_{2,\varepsilon}(g_{3,\varepsilon}(\delta)),\delta))) \cap \varp \subseteq
 \varepsilon^{\up}_{4,1}$ hence by clause (b)
\sn
\item[$(\beta)$]  if $\beta \in 
t_{g_{2,\varepsilon}(g_{3,\varepsilon}(\delta))}$ then
$\Rang(F_1(\rho_{\bar h}(\beta,\delta)) \cap
 \varepsilon \subseteq \varepsilon^{\up}_{4,1}$
\sn
\item[$(\gamma)$]   also $\Rang(F_0(\rho_{\bar h}(g_{2,\varepsilon}
(g_{3,\varepsilon}(\delta)),\delta))) \subseteq
 \varepsilon^{\up}_{4,0}$.
\end{enumerate}
\end{enumerate}
\end{enumerate}
\mn
[Why?  Recall that $\kappa_1$ is regular uncountable (being
$\theta_1$) and $\kappa_0 < \kappa_1$ is regular (being $\theta_0$).
First, for each $\delta \in \cU^{\up}_3$ we use Fodor lemma on
$S^{\kappa_1}_{\kappa_0,j_*} \backslash \varepsilon^{\up}_{1,1}$ to choose
$s_\delta,\varepsilon^{\up}_{4,1,\delta},\varepsilon^{\up}_{4,0,\delta}$
  as in clauses (c) + (d); second use the Fodor Lemma on $\cU^{\up}_3$ to get
 $\cU^{\up}_4,\varepsilon^{\up}_{4,1},\varepsilon^{\up}_{4,0}$; we
 cannot do it for $s_\delta$ as maybe $2^{\kappa_1} \ge \lambda$.

Let us verify $(d)(\beta)$ and $(d)(\gamma)$.  For $(d)(\beta)$ notice
that $\Rang(F_1(\rho_{\bar h}(\beta,\delta))) \subseteq
\varp^{\up}_{1,1} < \varp^{\up}_{4,1}$ for every 
$\beta \in t_{g_{2,\varp}(g_{3,\varp}(\delta))}$ by $(*)_1(d)$.  This
requirement is easy since $|t_{g_{2,\varp}(g_{3,\varp}(\delta))}| <
\kappa_1$ and $\rho_{\bar h}(\beta,\delta)$ is finite for every $\beta
\in t_{g_{2,\varp}(g_{3,\varp}(\delta))}$.

For $(d)(\gamma)$ we apply Fodor's lemma twice.

First, fix an ordinal $\delta \in \cU^{\up}_4$.  For every $\varp \in
s_\delta$, the sequence 
$F_0(\rho_{\bar h}(g_{2,\varp}(g_{3,\varp}(\delta)))$ is 
finite and hence bounded in
$\kappa_0$.  But $\kappa_0 < \kappa_1 = \cf(\kappa_1)$ and hence by
shrinking $s_\delta$ if needed we may assume that all the values are
bounded by the same ordinal $\sigma_\delta < \kappa_0$.

Now for each $\delta \in \cU^{\up}_4$ we choose $\sigma_\delta \in
\kappa_0$ in this way, so by shrinking $\cU^{\up}_4$ if needed we may
assume that $\sigma_\delta = \sigma$ for some fixed $\sigma <
\kappa_0$ and every $\delta \in \cU^{\up}_4$.  Now choose
$\varp^{\up}_{4,0} > \max\{\sigma,\varp^{\up}_{1,0}\}$.
\mn
\begin{enumerate}
\item[$(*)_5$]  we can find
  $\cU^{\dn}_1,\varepsilon^{\dn}_{1,0},\varepsilon^{\dn}_{1,1}$ such
  that:
\sn
\begin{enumerate}
\item[$(a)$]  $\cU^{\dn}_1 \subseteq S^*_0$ is stationary in $\lambda$
\sn
\item[$(b)$]  $\alpha < \delta \in \cU^{\dn}_1 \Rightarrow t_\alpha
  \subseteq \delta$
\sn
\item[$(c)$]  $\varepsilon^{\dn}_{1,\iota} < \kappa_\iota$ for $\iota
  = 0,1$
\sn
\item[$(d)$]  if $\delta \in \cU^{\dn}_1$ \then \, for arbitrarily large
$\alpha < \delta$ we have $\beta \in t^\iota_\alpha \wedge \iota \in
  \{0,1\} \Rightarrow \Rang(F_\iota(\rho_{\bar h}(\delta,\beta)))
\subseteq \varepsilon^{\dn}_{1,\iota} < \kappa_\iota$.
\end{enumerate}
\end{enumerate}
\mn
[Why?  Clearly $E = \{\delta < \lambda:\delta$ a limit ordinal such
  that $\alpha < \delta \Rightarrow t_\alpha \subseteq \delta\}$ is a
  club of $\lambda$.
For every $\delta \in S^*_0 \cap E$ and $\alpha < \delta$ we can find
$(\varepsilon^{\dn}_{1,0,\delta,\alpha},\varepsilon^{\dn}_{1,1,\delta,\alpha})$
as in clauses (c),(d) above because $|t^\iota_\alpha| < \kappa_\iota =
\cf(\kappa_\iota)$.  So recalling that 
$\cf(\delta) = \partial > \theta_1 = \kappa_1 > \kappa_0 = 
\theta_0$ it follows that there is a pair
$(\varepsilon^{\dn}_{1,0,\delta},\varepsilon^{\dn}_{1,1,\delta})$
  such that $\delta = \sup\{\alpha < \delta:
(\varepsilon^{\dn}_{1,0,\delta,\alpha},\varepsilon^{\dn}_{1,1,\delta,\alpha})
= (\varepsilon^{\dn}_{1,0,\delta},\varepsilon^{\dn}_{1,1,\delta})\}$.  
Then recalling $\lambda = \cf(\lambda) > \kappa_1 + \kappa_0$ we can
 choose $(\varepsilon^{\dn}_{1,0},\varepsilon^{\dn}_{1,1})$
 such that the set $\cU^{\dn}_1 = \{\delta \in S^*_0:
(\varepsilon^{\dn}_{1,0,\delta},\varepsilon^{\dn}_{1,1,\delta})
= (\varepsilon^{\dn}_{1,0},\varepsilon^{\dn}_{1,1})\}$ is stationary.]
\mn
\begin{enumerate}
\item[$(*)_6$]  we can find $\cU^{\dn}_2,\varepsilon^{\dn}_{2,0},
\varepsilon^{\dn}_{2,1}$ such that:
\sn
\begin{enumerate}
\item[$(a)$]  $\cU^{\dn}_2 \subseteq S^*_0 \backslash (\alpha^*_3 +1)$
  is stationary
\sn
\item[$(b)$]  if $\delta \in \cU^{\dn}_2$ and $\zeta < \kappa_1$
 \then \, $\delta = \sup(\cU^{\dn}_1 \cap \delta)$ and
for arbitrarily large $\delta_0 \in \cU^{\dn}_1 \cap \delta$ we have 
$\zeta < \max \Rang(F_1(\rho_{\bar h}(\delta,\delta_0)))$ 
and $\Rang(F_0(\rho_{\bar h}(\delta,\delta_0))) \subseteq 
\varepsilon^{\dn}_{2,0}$ and $\zeta \cap \Rang(F_1(\rho_{\bar
 h}(\delta,\delta_0))) \subseteq \varepsilon^{\dn}_{2,1}$
\sn
\item[$(c)$]   $\varepsilon^{\dn}_{2,0} \in
 (\varepsilon^{\dn}_{1,0},\kappa_0)$ and $\varepsilon^{\dn}_{2,1} \in
 (\varepsilon^{\dn}_{1,1},\kappa_1)$.
\end{enumerate}
\end{enumerate}
\mn
[Why?  For every $\zeta < \kappa_1$ let $S'_\zeta = \{\alpha \in S:\alpha =
  \sup(\cU^{\dn}_1 \cap \alpha)$ and $F_1(h(\alpha)) = \zeta\}$,
  clearly it is a stationary subset of $\lambda$.

Let $E = \{\delta < \lambda:\delta$ is a limit ordinal and $\zeta <
\kappa_1 \Rightarrow \delta = \sup(\delta \cap S'_\zeta)\}$. 
Clearly it is a club of $\lambda$.
If $\zeta \in S^{\kappa_1}_{\kappa_0} \backslash \varp^{\dn}_{1,1}$ and
$\delta \in E \cap S^*_0$ and $\alpha \in S'_\zeta \cap \delta$ let
 $\varepsilon^{\dn}_{2,0,\zeta,\delta,\alpha} =
 \sup\Rang(F_0(\rho_{\bar h}(\delta,\alpha))) + \varp^{\dn}_{1,0} + 1$ and let 
$\varepsilon^{\dn}_{2,1,\zeta,\delta,\alpha} = \sup(\zeta \cap
 \Rang(F_1(\rho_{\bar h}(\delta,\alpha))) +1 < \zeta$.

Fixing $\delta$ and $\zeta$, recalling $\cf(\delta) > \kappa_0 + \kappa_1$, 
for some pair $(\varepsilon^{\dn}_{2,0,\zeta,\delta},
\varepsilon^{\dn}_{2,1,\zeta,\delta}) \in \kappa_0 \times \kappa_1$
we have $\delta = \sup\{\alpha \in S'_\zeta \cap \delta:
(\varepsilon^{\dn}_{2,0,\zeta,\delta,\alpha},
\varepsilon^{\dn}_{2,1,\zeta,\delta,\alpha}) =
(\varepsilon^{\dn}_{2,0,\zeta,\delta},
\varepsilon^{\dn}_{2,1,\zeta,\delta})\}$.

Fixing $\delta$ apply Fodor lemma on $S^{\kappa_1}_{\kappa_0}$, for
some pair $(\varepsilon^{\dn}_{2,0,\delta},\varepsilon^{\dn}_{2,1,\delta})$ the
set $b_\delta = \{\zeta \in S^{\kappa_1}_{\kappa_0}:
(\varepsilon^{\dn}_{2,0,\zeta,\delta},
\varepsilon^{\dn}_{2,1,\zeta,\delta}) =
(\varepsilon^{\dn}_{2,0,\delta},\varepsilon^{\dn}_{2,1,\delta})\}$ is
a stationary subset of $\kappa_1$.  

Applying Fodor lemma on $\delta \in E \cap S^*_0$, there is a pair
$(\varepsilon^{\dn}_{2,0},\varepsilon^{\dn}_{2,1})$ such that
$\cU^{\dn}_2 := \{\delta \in S^*_0:\delta \in E$ and
$(\varepsilon^{\dn}_{2,0,\delta},\varepsilon^{\dn}_{2,1,\delta}) =
(\varepsilon^{\dn}_{2,0},\varepsilon^{\dn}_{2,1})\}$ is stationary.
Clearly we are done.  We could have put $b_\varp$ in $(*)_6(b)$ but it
does not seem needed.]
\medskip

\noindent
\underline{Stage C}:  Now we shall find the required $\alpha_0 < \alpha_1$.

In this stage we deal with $\bold c_1$, so $j_* < \kappa_1$.
First, there are
$\delta_1,\delta_2,\varepsilon^{\md}_0,\varepsilon^{\md}_1,\alpha^*_4$
such that:
\mn
\begin{enumerate}
\item[$\oplus_0$]  
\begin{enumerate}
\item[(a)]  $\delta_1 \in \cU^{\dn}_2$ and $\delta_2 \in \cU^{\up}_4$,
  see $(*)_6$ and $(*)_4$ respectively
\sn
\item[(b)]  $\delta_1 < \delta_2$ and $\alpha^*_3 < \delta_1$
\sn
\item[(c)]  $\varepsilon^{\md}_\iota :=
\max\Rang(F_\iota(\rho_{\bar h}(\delta_2,\delta_1))) > 
\varepsilon^{\dn}_{2,\iota} + \varepsilon^{\up}_{4,\iota} \ge
\varp^{\dn}_{1,\iota} + \varp^{\up}_{1,\iota}$ for $\iota=0,1$
\sn
\item[(d)]  $\alpha^*_4 < \delta_1$ is $> \alpha^*_3$ 
and if $\alpha \in (\alpha^*_4,\delta_1)$ 
then $\rho_{\delta_2,\delta_1} \char 94  
\langle \delta_1 \rangle \trianglelefteq \rho_{\delta_2,\alpha}$.
\end{enumerate}
\end{enumerate}
\mn
[Why can we?  Easy but we give details.  First, let $\cW_* = \{\delta
\in S:\delta$ is a limit ordinal $> \alpha^*_3$ necessarily of cofinality
  $\partial$ such that $F_\iota(\delta) > \varepsilon^{\dn}_{2,\iota} +
  \varepsilon^{\up}_{4,\iota}$ for $\iota = 0,1$ and $\delta = 
\sup(\delta \cap \cU^{\dn}_2)\}$, 
clearly it is a stationary subset of $\lambda$.
  Second, choose $\delta_2 \in \cU^{\up}_4$ which is $> \alpha^*_3$
  such that $\delta_2 = \sup(\cW_* \cap \delta_2)$.  Third, choose
$\delta_* \in \cW_* \cap \delta_2$ such that $\alpha^*_3 < \delta_*$.  
Fourth, let $\alpha_* < \delta_*$ be such that $\alpha_* > \alpha^*_3$ and 
$\alpha \in (\alpha_*,\delta_*) \Rightarrow \rho(\delta_2,\delta_*) 
\char 94 \langle \delta_* \rangle 
\trianglelefteq \rho(\delta_2,\alpha)$ (hence $\rho_{\bar h}
(\delta_2,\delta_*) \char 94 \langle h_{\delta_* + 1}(\delta_*)\rangle 
\trianglelefteq \rho_{\bar h}(\delta_2,\alpha))$.
Fifth, choose $\delta_1 \in
 (\alpha_*,\delta_*) \cap \cU^{\dn}_2$ hence $\delta_1 > \alpha^*_3$.  
Sixth, we choose
  $\varepsilon^{\md}_\iota$ for $\iota = 0,1$ by clause (c), the
  inequality holds because $\delta_* \in \cW_* 
\cap \Rang(\rho_{\bar h}(\delta_2,\delta_1))$.  

Lastly, choose $\alpha^*_4$ as in $\oplus_0(d)$. 
Easy to check that we are done proving $\oplus_0$.]

Let $\rho = \rho_{\bar h}(\delta_2,\delta_1)$. 

Second, choose $\delta_0$ such that
\mn
\begin{enumerate}
\item[$\oplus_{0.1}$] 
\begin{enumerate}
\item[(a)]  $\delta_0 \in \cU^{\dn}_1 \cap \delta_1$
\sn
\item[(b)]  $(*)_6(b)$ holds with
$(\varepsilon^{\md}_1,\delta_1)$ here standing for $(\zeta,\delta)$
there, that is, we have $\varp^{\md}_1 < \max \Rang(F_1(\rho_{\bar
  h}(\delta_1,\delta_0)))$ and $\Rang(F_0(\rho_{\bar
  h}(\delta_1,\delta_0))) \subseteq \varp^{\dn}_{2,0}$ and
$\varp^{\md}_1 \cap \Rang(F_1(\rho_{\bar h}(\delta_1,\delta_0)))
\subseteq \varp^{\dn}_{2,1}$
\sn
\item[(c)]  $\delta_0 > \alpha^*_4$ recalling $\delta_1 >
  \alpha^*_4 > \alpha^*_3$ by $\oplus_0(b),(d)$.
\end{enumerate}
\end{enumerate}
\mn
[Why can we choose $\delta_0$?  By $(*)_6$.]

Also choose $\alpha^*_5$ such that
\mn
\begin{enumerate}
\item[$\oplus_{0.2}$]  $\alpha^*_5 < \delta_0$ is such that $\alpha \in
 (\alpha^*_5,\delta_0) \Rightarrow \rho_{\delta_1,\delta_0}
 \char 94 \langle \delta_0 \rangle \trianglelefteq
 \rho_{\delta_1,\alpha}$.
\end{enumerate}
\mn
Third, choose $\varepsilon_* \in s_{\delta_2}$ ($s_{\delta_2}$ is from
 $(*)_4(c),(d))$ such that
$\varepsilon_* > \varepsilon^{\md}_{2,1} := \max \big(\Rang(F_1(\rho_{\bar
  h}(\delta_2,\delta_1) \cup \Rang(F_1(\rho_{\bar
  h}(\delta_1,\delta_0)))\big)$ which is $> \varepsilon^{\md}_1$, 
possible as $s_{\delta_2}$ is a stationary subset of $\kappa_1$.  

Fourth, let $\delta_3 = g_{3,\varepsilon_*}(\delta_2)$.  

Fifth, let $\alpha_1 = g_{2,\varepsilon_*}(\delta_3)$.

Lastly, choose $\alpha_0 < \delta_0$ large enough and as in $(*)_5(d)$
such that $\alpha_0 > \alpha^*_5 > \alpha^*_4$, that is, we have
$\beta \in t^1_{\alpha_0} \wedge 1 \in \{0,1\} \Rightarrow
\Rang(F_1(\rho_{\bar h}(\delta_0,\beta))) \subseteq
\varp^{\dn}_{1,1} < \kappa_1$, used only for $1 = 0$.

We shall prove below that the pair $(\alpha_0,\alpha_1)$ is as
promised.  

So (finishing the case of $\kappa_1$ colours)
\mn
\begin{enumerate}
\item[$\circledast$]  let $\zeta_0 \in 
t^0_{\alpha_0},\zeta_1 \in t^1_{\alpha_1}$ 
and we should prove that $\bold c_1\{\zeta_0,\zeta_1\} = j_*$.
\end{enumerate}
\mn
Note
\mn
\begin{enumerate}
\item[$\oplus_1$]  $\delta_2 \in \cU^{\up}_4 \subseteq \cU^{\up}_3$
  and $\alpha_0 < \delta_0 < \delta_1 < \delta_2$.
\end{enumerate}
\mn
[Why?  The first statement holds by 
the choice of $\delta_2$, see $\oplus_0(a)$ and $(*)_4(a)$.  The
second statement holds by the choices of $\delta_1$,
i.e. $\oplus_0(b)$, the choice of $\delta_0$, i.e. $\oplus_{0.1}(a)$
and the choice of $\alpha_0$ (see ``Lastly..." after $\oplus_{0.2}$).]
\mn
\begin{enumerate}
\item[$\oplus_2$]  $\delta_3 = g_{3,\varepsilon_*}(\delta_2) \in 
\cU^{\up}_{2,\varepsilon_*}$ and $\delta_2 < \delta_3$.
\end{enumerate}
\mn
[Why?  By the choice of $\delta_3$ (after $\oplus_{0.2}$ in ``Fourth") and by 
$(*)_3(d)+(e)$ (note that the assumption of $(*)_3(e)$ in our case,
which means $\delta_2 \in \cU^{\up}_3$ and $\varp_* \in
S^{\kappa_1}_{\kappa_0} \backslash \varp^{\md}_{2,1}$, holds by $\oplus_1$
and by the ``Third" after $\oplus_{0.2}$ above (recalling
$s_{\delta_2} \subseteq S^{\kappa_1}_{\kappa_0}$ and $\oplus_0(c)$)).]
\mn
\begin{enumerate}
\item[$\oplus_3$]  $\alpha_1 = g_{2,\varepsilon_*}(\delta_3) 
\in \cU^{\up}_1$ and $\delta_3 < \alpha_1$.
\end{enumerate}
\mn
[Why?  By the choice of $\alpha_1$ in ``Fifth" after $\oplus_{0,2}$ 
and $(*)_2(d)$.]
\mn
\begin{enumerate}
\item[$\oplus_4$]  $\eta_0 := \rho_{\bar h}(\zeta_1,\alpha_1)$
satisfies ($\eta_0 \in {}^{\omega >}\lambda$ and):
\sn
\begin{itemize}
\item   $\Rang(F_1(\eta_0)) \subseteq \varp^{\up}_{1,1}$.
\end{itemize}
\end{enumerate}
\mn
[Why?  By $(*)_1(d)$ recalling $\oplus_3$ and, of course,
$\alpha_1 > \alpha^*_5 > \alpha^*_1$ the case $\iota=0$ is 
unreasonable as $t^1_{\alpha_1}$ may be of cardinality $\ge \kappa_0 
= \theta$.]

Recall that $(*)_1(d)$ deals only with $t^1_\varp$.
\mn
\begin{enumerate}
\item[$\oplus_5$]  $\nu_0 := \rho_{\bar h}(\alpha_1,\delta_2)$
  satisfies ($\nu_0 \in {}^{\omega >}\lambda$ and)
\sn
\begin{enumerate}
\item[(a)]  $\Rang(F_0(\nu_0)) \subseteq \varepsilon^{\up}_{4,0}$
\sn
\item[(b)]  $\varepsilon_* \in \Rang(F_1(\nu_0))$
\sn
\item[(c)]  $\Rang(F_1(\nu_0)) \cap \varepsilon_*
  \subseteq \varepsilon^{\up}_{4,1}$
\sn
\item[(d)]  $\alpha_1 =
  g_{2,\varepsilon_*}(g_{3,\varepsilon_*}(\delta_2)) =
  g_{2,\varp_*}(\delta_3)$ 
\sn
\item[(e)]  $\rho(\alpha_1,\delta_2) =
  \rho(\alpha_1,g_{3,\varp_*}(\delta_2)) \char 94 \rho 
(g_{3,\varepsilon_*}(\delta_2),\delta_2)$.
\end{enumerate}
\end{enumerate}
\mn
[Why?  Clause (d) of $\oplus_5$ holds by the choice of $\alpha_1$ in
``Fourth" and ``Fifth" after $\oplus_{0.2}$ above (and see $\oplus_2$); 
similarly clause (e) holds.
By $\oplus_1$ we have $\delta_2 \in \cU^{\up}_4$ and by
$(*)_4(d)(\gamma),(\alpha)$ and the choices of 
$\delta_3,\alpha_1$ we have clauses (a) + (c) 
of $\oplus_5$; that is, $(\alpha_1,\delta_{2,\varp_*},\varp_*)$ here stand for
$(g_{2,\varp}(g_{3,\varp}(\delta)),\delta,\varp)$ in $(*)_4(d)$.  
Now $\delta_3 \in \Rang(\rho(g_{3,\varp_*}(\delta)),\delta_2)$ by
$\oplus_2$ hence $\delta_3 \in \Rang(\rho(\alpha_1,\delta_2))$ by
$\oplus_5(e)$ hence $\delta_3 \in \Rang(\nu_0)$ by the choice of
$\nu_0$ (see the beginning of $\oplus_5$).  This implies clause (b) of
$\oplus_5$ because $F_1(\delta_3) = \varp_*$ because $\delta_3 \in
\dom(g_{2,\varp_*}) \subseteq \cU^{\up}_{2,\varp_*}$ by $\oplus_2$ and
$(\forall \delta)[\delta \in \cU^{\up}_{2,\varp_*} \Rightarrow \delta
\in S^*_{\gamma^*_{\varp_*}} 
\Rightarrow F_1(\delta) = \varp_*]$ by $(*)_2(a),(b)$.]
\mn
\begin{enumerate}
\item[$\oplus_6$]  $\nu_1 := \rho_{\bar h}(\delta_1,\delta_0)$ satisfies:
\sn
\item[${{}}$]  $(a) \quad \Rang(F_0(\nu_1)) \subseteq \varepsilon^{\dn}_{2,0}$
\sn
\item[${{}}$]  $(b) \quad \varepsilon^{\md}_1 < \max \Rang(F_1(\nu_1))$
\sn
\item[${{}}$]  $(c) \quad \Rang(F_1(\nu_1)) \subseteq \varepsilon_*$.
\end{enumerate}
\mn
[Why?  By $\oplus_0(a)$ we have $\delta_1 \in \cU^{\dn}_2$.  So (a),(b) 
hold by $(*)_6(b)$ and the choice of $\delta_0$, i.e. $\oplus_{0.1}(b)$;
 we use the
first two conclusions of $(*)_6(b)$ not the third.  As for clause (c) it
holds by the choice of $\varepsilon_*$ in ``Third" after $\oplus_{0.2}$.]
\mn
\begin{enumerate}
\item[$\oplus_7$]  
\begin{enumerate}
\item[(a)]  $\eta_1 := \rho_{\bar h}(\delta_0,\zeta_0)$ satisfies
\sn
\begin{itemize}
\item  $\Rang(F_\iota(\eta_1)) \subseteq
  \varepsilon^{\dn}_{1,\iota}$ for $\iota = 0,1$.
\end{itemize}
\sn
\item[(b)]  $\rho = \rho_{\bar h}(\delta_2,\delta_1)$ satisfies
\sn
\begin{itemize}
\item  $\max \Rang(F_\iota(\rho)) = \varepsilon^{\md}_\iota$ for $\iota=0,1$.
\end{itemize}
\end{enumerate}
\end{enumerate}
\mn
[Why?  Clause (a) holds by $(*)_5(d)$ and the choice of 
$\alpha_0$ in ``lastly" after
$\oplus_{0.2}$ recalling $\zeta_0 \in t^0_{\alpha_0}$.  Clause (b)
holds by $\oplus_0(c)$.]
\mn
\begin{enumerate}
\item[$\oplus_8$]
\begin{enumerate}
\item[(a)]  $\rho_{\bar h}(\zeta_1,\zeta_0) =
\rho_{\bar h}(\zeta_1,\alpha_1) \char 94 \rho_{\bar h}(\alpha_1,\delta_2) 
\char 94 \rho_{\bar h}(\delta_2,\delta_1) \char 94 
\rho_{\bar h}(\delta_1,\delta_0) \char 94 \rho_{\bar h}(\delta_0,\zeta_0)$
\sn
\item[(b)]  recalling $\rho = \rho_{\bar h}
(\delta_2,\delta_1)$ and the choices of $\eta_0,\nu_0,\rho,\nu_1,\eta_1$
we have $\rho_{\bar h}(\zeta_1,\zeta_0) = \eta_0 \char 94
  \nu_0 \char 94 \rho \char 94 \nu_1 \char 94 \eta_1$.
\end{enumerate}
\end{enumerate}
\mn
[Why?  Clause (a) holsd by the choices of $\alpha^*_0$ in $(*)_1(c)(d)$ and
  of $\alpha^*_3$ in $(*)_3(f),(f)'$ and $\delta_1 > \alpha^*_3$ by
  $\oplus_0(b)$ and as ``$\delta_0 > \alpha^*_3$" recalling
  $\oplus_{0.1}(c)$ and ``$\alpha_0 > \alpha^*_5$", see ``Lastly"
 after $\oplus_{0.2}$.  Clause (b)
  holds by clause (a) and the definitions of
  $\eta_0,\nu_0,\rho,\nu_1,\eta_1$ above, that is, in $\oplus_4$, in
  $\oplus_3$, before $\oplus_{0.1}$, in $\oplus_6$, in $\oplus_7$
  respectively.]
\mn
\begin{enumerate}
\item[$\oplus_9$]  $\ell^\bullet_4 := \bold d(\rho_{\bar h}(\zeta_1,\zeta_0))$
 satisfies $F_1(\varrho(\ell^\bullet_4)) = \varepsilon_*$.
\end{enumerate}
\mn
[Why?  We shall use $\oplus_8(a),(b)$ freely; now $\bold d$ was chosen by Claim
\ref{e4} and letting $\varrho = \eta_0 \char 94 \nu_0 \char 94 \rho
  \char 94 \nu_1 \char 94 \eta_1$ we apply the claim to
  $(\eta_0,\nu_0,\rho,\nu_1,\eta_1)$, so it suffices to show that
clauses (B)(a)-(d) of \ref{e4} hold.
\mn
\begin{enumerate}
\item[$\oplus_{9.1}$]  clause $(B)(a)(\alpha)$ of \ref{e4} holds.
\end{enumerate}
\mn
Why?  First, $\varepsilon_* \le \max \Rang(F_1(\nu_0))$ by
$\oplus_5(b)$.  

Second, $\Rang(F_1(\eta_0)) \subseteq
\varepsilon^{\up}_{1,1}$ by $\oplus_4$ and $\varepsilon^{\up}_{1,1}
\le \varepsilon^{\up}_{4,1}$ by $(*)_4(b)$ and
$\varepsilon^{\up}_{4,1} \le \varepsilon^{\md}_1$ by $\oplus_0(c)$ and
$\varepsilon^{\md}_1 < \varepsilon_*$ by the choice of $\varepsilon_*$
in ``Third" after $\oplus_{0.2}$.  

Third, $\Rang(F_1(\rho)) \subseteq \varepsilon_*$
as $\Rang(F_1(\rho)) = \Rang(F_1(\rho_{\bar h}(\delta_2,\delta_1))) 
\subseteq \varepsilon^{\md}_1 +1$ by $\oplus_0(c)$ and
$\varepsilon^{\md}_1 < \varepsilon_*$ by the choice of $\varepsilon_*$.  

Fourth, $\Rang(F_1(\nu_1)) \subseteq \varepsilon_*$ by $\oplus_6(c)$.  

Fifth, $\Rang(F_1(\eta_1)) \subseteq \varepsilon_*$
as $\Rang(F_1(\eta_1)) \subseteq \varepsilon^{\dn}_{1,1}$ by
$(*)_5$ and $\varepsilon^{\dn}_{1,1}
\le \varepsilon^{\dn}_{2,1}$ by $(*)_6(c)$ and
$\varepsilon^{\dn}_{2,1} < \varepsilon^{\md}_1$ by $\oplus_0(c)$ and
$\varepsilon^{\md}_1 < \varepsilon_*$ by the choice of
$\varepsilon_*$.

Together $\oplus_{9.1}$ holds.
\mn
\begin{enumerate}
\item[$\oplus_{9.2}$]  let $\ell_1 < \ell g(\nu_0)$ 
be as in clause $(B)(a)(\beta)$ of \ref{e4}
\sn
\item[$\oplus_{9.3}$]  clause $(B)(b)(\alpha)$ of \ref{e4} holds.
\end{enumerate}
\mn
Why?  First, $\max \Rang(F_0(\rho)) = \varepsilon^{\md}_0$ by $\oplus_0(c)$.  

Second, $\Rang(F_0(\eta_0)) \subseteq \varepsilon^{\md}_0$ is
unreasonable see $\oplus_4$ and not necessary.

Third, $\Rang(F_0(\nu_0)) \subseteq \varepsilon^{\md}_0$ because
$\Rang(F_0(\nu_0)) \subseteq \varepsilon^{\up}_{4,0}$ by $\oplus_5(a)$ and
$\varepsilon^{\up}_{4,0} \le \varepsilon^{\md}_0$ by $\oplus_0(c)$.

Fourth, $\Rang(F_0(\nu_1)) \subseteq \varepsilon^{\md}_0$ because
$\Rang(F_0(\nu_1)) \subseteq \varepsilon^{\dn}_{2,0}$ by
$\oplus_6(a)$ and $\varepsilon^{\dn}_{2,0} \le \varepsilon^{\md}_0$ by
$\oplus_0(c)$.

Fifth, $\Rang(F_0(\eta_1)) \subseteq \varepsilon^{\md}_0$
because $\Rang(F_0(\eta_1)) \subseteq \varepsilon^{\dn}_{1,0}$ by
$\oplus_7(a)$ and $\varepsilon^{\dn}_{1,0} < \varepsilon^{\dn}_{2,0}$ by
$(*)_6(c)$ and $\varepsilon^{\dn}_{2,0} \le \varepsilon^{\md}_0$ by
$\oplus_0(c)$.

Together $\oplus_{9.3}$ holds.
\mn
\begin{enumerate}
\item[$\oplus_{9.4}$] 
\begin{enumerate}
\item[(a)]   let $\ell^\bullet_2 < \ell g(\varrho)$ 
be as in clause $(B)(b)(\beta)$ of \ref{e4}
\sn
\item[(b)]   let $\ell^*_2 = \ell^\bullet_2 - \ell g(\eta_0 \char 94 \nu_0)$
\end{enumerate}
\sn
\item[$\oplus_{9.5}$] 
\begin{enumerate}
\item[(a)]   $\ell^\bullet_2  \in [\ell g(\eta_0 \char 94 \nu_0),\ell g(\eta_0
  \char 94 \nu_0 \char 94 \rho))$
\sn
\item[(b)]   clause $(B)(c)(\alpha)$ holds, i.e.
\sn
\begin{enumerate}
\item[$\bullet_1$]  $\max\Rang(F_1(\nu_0)) > \max\Rang(F_1(\varrho
  \rest [\ell^\bullet_2,\ell g(\varrho)))$
\sn
\item[$\bullet_2$]  $\max\Rang(F_1)(\varrho \rest [\ell^\bullet_2,\ell
  g(\varrho))) = \max\Rang(F_1(\nu_1)) > \max \Rang(\rho \char 94 \eta_1)$
\end{enumerate}
\sn
\item[(c)] let $\ell_3 < 
\ell g(\nu_1)$ be as in clause $(B)(c)(\beta)$ of \ref{e4} 
\sn
\item[(d)]  $F_1(\nu_1(\ell_3)) \ge \varp^{\md}_1$.
\end{enumerate}
\end{enumerate}
\mn
Why?  Clause (a) follows by $(B)(b)(\alpha)$ proved in $\oplus_{9.3}$ above.
Clause $(b),\bullet_1$ holds by $\oplus_{9.1}$.  Clause $(b),\bullet_2$
follows because: first $\Rang(F_1(\rho))
\subseteq \varp^1_{\md} +1$ by $\oplus_0(c)$ and $\varp^1_{\md} +1 <
\varp$ by second; $\Rang(F_1(\nu_1)) \nsubseteq \varp^{\md}_1+1$ by
$\oplus_6(b)$ and third, $\Rang(F_1(\eta_1)) \subseteq
\varp^{\dn}_{1,1}$ by $\oplus_7(a)$ and $\varp^{\dn}_{1.1} <
\varp^{\md}_1$ by $\oplus_0(d)$ by the choice of $\varp_*$.

By clause (b), it follows that $\ell_3$ from Clause (c)
are well defined and Clause (d) holds
\mn
\begin{enumerate}
\item[$\oplus_{9.6}$] 
\begin{enumerate}
\item[(a)]  $\Rang(F_1(\eta_0 \char 94 (\rho \rest \ell^*_2) \char 94
  \nu_1 \char 94 \eta_1)) \subseteq \varepsilon^{\md}_1 +1$
\sn
\item[(b)]  $\varepsilon_* \in \Rang(F_1(\nu_0))$ is $> \varp^{\md}_1$
\sn
\item[(c)]  $\Rang(F_1(\nu_0)) \cap \varepsilon_* 
\subseteq \varepsilon^{\md}_1$
\end{enumerate}
\end{enumerate}
\mn
Why?  First, $\Rang(F_1(\eta_0)) \subseteq \varepsilon^{\md}_1$
because $\Rang(F_1(\eta_0)) \subseteq \varepsilon^{\up}_{1,1}$ by
$\oplus_4$ and $\varepsilon^{\up}_{1,1} \le \varepsilon^{\up}_{4,1}$ by
$(*)_4(b)$ and $\varepsilon^{\up}_{4,1} \le \varepsilon^{\md}_1$ by
$\oplus_0(c)$.

Second, $\Rang(F_1(\rho \rest \ell^*_2)) \subseteq \Rang(\rho)
\subseteq \varp^{\md}_1 +1$ and $\Rang(F_1(\nu_1 \char 94 \eta_1))
\subseteq \varp^{\md}_1$ by $\oplus_0(c)$.

Third, $\Rang(F_1(\eta_0 \char 94 (\rho \rest \ell^*_2)) \subseteq 
\Rang(F_1(\eta_0)) \cup \Rang(F_1(\rho \rest \ell^*_2)) \subseteq 
\varepsilon^{\md}_1 +1$ by the last two sentences, so clause (a) of
$\oplus_{9.6}$ holds.

Fourth, clause (b), i.e. $\varepsilon_* \in \Rang(F_1(\nu_0))$ holds
by $\oplus_5(b)$.

Fifth, $\Rang(F_1(\eta_0 \char 94 \nu_0)) \cap \varepsilon_* \subseteq
\varepsilon^{\up}_{4,1}$ by $(*)_4(d)$ with $(\delta,\beta,\varp)$ there
standing for $(\delta_2,\zeta_1,\varp_*)$ here (recalling $\delta_2 \in
\cU^{\up}_4$ and $\zeta_1 \in t^1_{\alpha_1} =
t^1_{g_{2,\varepsilon_*}(g_{3,\varepsilon_*}(\delta_2))})$ and
$\varepsilon^{\up}_{4,1} \le \varepsilon^{\md}_1$ by $\oplus_0(c)$.
Hence, $\Rang(F_1(\nu_0)) \cap \varepsilon_* \subseteq
\Rang(F_1(\eta_0 \char 94 \nu_0)) \cap \varepsilon_* \subseteq
\varepsilon^{\up}_{4,1} \subseteq \varepsilon^{\md}_1$, so also clause
(c) of $\oplus_{9.6}$ holds. 
\mn
\begin{enumerate}
\item[$\oplus_{9.7}$]  
\begin{enumerate}
\item[(a)]   let $\ell^\bullet_4$ from $\oplus_9$ be as in $(B)(d)(\beta)$
\sn
\item[(b)]  $F_1(\varrho(\ell^\bullet_4)) = \varepsilon_*$
\sn
\item[(c)]  (used in stage D) $\ell^\bullet_4 \in [\ell g(\eta_0),\ell
  g(\eta_0 \char 94 \nu_1))$.
\end{enumerate}
\end{enumerate}
\mn
[Why?  By $\oplus_{9.6},\, \ell^\bullet_4$ is well defined and belongs
to $[\ell g(\eta_0),\ell g(\eta_0 \char 94 \nu_0))$; moreover,
$F_1(\varrho(\ell^\bullet_4)) = \varp_*$.]

So indeed $\oplus_9$ holds.
\mn
\begin{enumerate}
\item[$\oplus_{10}$]  $\bold c_1\{\zeta_0,\zeta_1\} = j_*$.
\end{enumerate}
\mn
[Why?  Because $\bold d(\varrho) = \ell^\bullet_4$ and
 $(F_1(\varrho))(\ell^\bullet_4) = \varepsilon_*$ and so by
  $\odot_7(c), h''(\varp_*) = \ell^\bullet_4$ 
we have $\bold c_1\{\zeta_0,\zeta_1\} =
 h'(\varepsilon_*)$ and $h'(\varepsilon_*) =j_*$ because
 $\varepsilon_* \in s_{\delta_2}$ by the choice of $\varepsilon_*$
and $h'(\varepsilon_*)$ is $j_*$ by $(*)_4(c)$ recalling the
definition of $S^{\kappa_1}_{\kappa_0,j_*}$ in $\odot_7(a)$.]
\bigskip

\noindent
\underline{Stage D}:

We would like to have $\lambda$ colours (not just $\kappa_1$ colours),
but (unlike earlier versions) we rely on what was proved (i.e. the
properties of $\bold c_1$) instead of repeating it. 
So we shall assume $\boxplus$ from the beginning of Stage B and $j_* <
\lambda$ in $\boxplus(d)$.

Now
\mn
\begin{enumerate}
\item[$\boxplus_1$]  for some $\cW_1,\varepsilon^{\up}_{0,1},\alpha^*_{0,1}$ 
\sn
\begin{enumerate}
\item[(a)]  $\alpha^*_{0,1} < \lambda,\varepsilon^{\up}_{0,1} < \kappa_1$
\sn
\item[(b)]  $\cW_1 \subseteq S$ is stationary and  
$\min(\cW_1) > \alpha^*_{0,1}$
\sn
\item[(c)]   if $\delta \in \cW_1$ and $\beta \in t_\delta$
  then $\Rang(F_1(\rho_{\bar h}(\beta,\delta))) \subseteq
  \varepsilon^{\up}_{0,1}$
\sn
\item[(d)]  if $\delta \in \cW_1$ and $\alpha \in 
[\alpha^*_{0,1},\delta)$ and $\beta \in t_\delta$ then
 $\rho(\beta,\delta) \char 94 \langle \delta \rangle \trianglelefteq
 \rho(\beta,\alpha)$.
\end{enumerate}
\end{enumerate}
\mn
[Why?  As in the proof of $(*)_1$ in Stage B.]
\mn
\begin{enumerate}
\item[$\boxplus_2$]
\begin{enumerate}
\item[(a)]  let $\cW_2 = \{\delta \in S:F_2(h(\delta)) 
= j_*,F_1(h(\delta)) = \varepsilon^{\up}_{0,1}$ and
 $\delta > \alpha^*_{0,1}\}$, so stationary
\sn
\item[(b)]   let $g^*_1:\cW_2 \rightarrow \cW_1$ be such
 that $\delta < g^*_1(\delta) \in \cW_1$
\end{enumerate}
\sn
\item[$\boxplus_3$]  there are $\cW_3,\alpha^*_{0,2}$ and $n_*$ such that:
\sn
\begin{enumerate}
\item[(a)]  $\cW_3 \subseteq \cW_2$ is stationary and
  $\min(\cW_3) > \alpha^*_{0,2} > \alpha^*_{0,1}$
\sn
\item[(b)]   if $\delta \in \cW_3$ and $\alpha \in 
[\alpha^*_{0,2},\delta)$ and $\beta \in t_{g^*_1(\delta)}$ then
 $\rho(\beta,g^*_1(\delta)) \char 94 \langle g^*_1(\delta)\rangle 
\trianglelefteq \rho(\beta,g^*_1(\delta)) \char 94
 \rho(g^*_1(\delta),\delta) \char 94 \langle \delta\rangle
 \trianglelefteq \rho(\beta,\alpha)$
\sn
\item[(c)]   if $\delta \in \cW_3$ and $\beta \in t_{g^*_1(\delta)}$ \then \,
\sn
\begin{enumerate}
\item[$(\alpha)$]  $\Rang(F_1(\rho_{\bar h}(\beta,g^*_1(\delta))) 
\subseteq \varepsilon^{\up}_{0,1}$
\sn
\item[$(\beta)$]  $n_* = |\{\ell < k(\beta,\delta):
(F_1(\rho_{\bar h}(\beta,\delta))(\ell) = \varepsilon^{\up}_{0,1}\}|$
\sn
\item[$(\gamma)$]   hence if $\alpha < \delta$ and
 $\rho(\beta,\delta) \char 94 \langle \delta \rangle \trianglelefteq
 \rho(\beta,\alpha)$ then the $(n_* +1)$-th member of the set $\{\ell <
  k(\beta,\alpha):F_1(\rho_{\bar h}(\beta,\alpha))(\ell) =
\varepsilon^{\up}_{0,1}\}$ is $\ell g(\rho(\beta,\delta))$.
\end{enumerate}
\end{enumerate}
\end{enumerate}
\mn
[Why?  As usual, e.g. how do we justify $n_*$ in clause $(c)(\beta)$
not depending on $\beta \in t_\delta$?  First, find $\delta$, then for
any $\beta \in t_\delta$ we have
\mn
\begin{enumerate}
\item[$\bullet$]  $\rho(\beta,\delta) = \rho(\beta,g^*_1(\delta))
  \char 94 \rho(g^*_1(\delta),\delta)$.
\end{enumerate}
\mn
[Why?  Recall $\boxplus_1(d)$.]
\mn
\begin{enumerate}
\item[$\bullet$]  $\Rang(F_1(\rho_{\bar h}(\beta,g^*_1(\delta)))
  \subseteq \varepsilon^{\up}_{0,1}$.
\end{enumerate}
\mn
[Why?  Recall $\boxplus_1(c)$.]

Together, $n_*$ depends just on $\rho_{\bar h}(g^*_1(\delta),\delta)$
which does not depend on $\delta$.  Second, as choosing $\cW_3$ we can
make $n_*$ not depend on $\delta$.]

Let $j_{**} < \kappa_1$ be such that $h'_1(j_{**}) =
\varepsilon^{\up}_{0,1},h'_2(j_{**}) = n_*$.   Next let $g_*:\lambda
\rightarrow \cW_3$ be increasing and define $s_\alpha =
t_{g_*(\alpha)},s^\iota_\alpha = t^\iota_{g_*(\alpha)}$ for $\iota=0,1$.
Now by what was proved in the earlier stages we can find $\alpha_0 <
\alpha_1 < \lambda$ such that if $\zeta_0 \in s^0_{\alpha_0} \wedge
\zeta_1 \in s^1_{\alpha_1}$ then $\bold c_1\{\zeta_0,\zeta_1\} =
j_{**}$.

Let $(\zeta_0,\zeta_1) \in s^0_{\alpha_0} \times s^1_{\alpha_1}$.  By the
choice of $\bold c_1$, in $\odot_7$ we have $\bold c_2$ from $\odot_9$
and by $\boxplus_3(c)(\gamma)$ we have $\bold c_2(\{\zeta_0,\zeta_1\}) =
j_*$.  But $(s^0_{\alpha_0},s^1_{\alpha_1}) =
(t^0_{g_*(\alpha_0)},t^1_{g_*(\alpha_1)})$ so $\alpha'_0 =
g_*(\alpha_0),\alpha'_1 = g_*(\alpha_1)$ are as required.
\end{PROOF}
\newpage

\bibliographystyle{alphacolon}
\bibliography{lista,listb,listc,listd,liste,listf,listv,listx,listy,listz}

\def\germ{\frak} \def\scr{\cal} \ifx\documentclass\undefinedcs
  \def\bf{\fam\bffam\tenbf}\def\rm{\fam0\tenrm}\fi 
  \def\defaultdefine#1#2{\expandafter\ifx\csname#1\endcsname\relax
  \expandafter\def\csname#1\endcsname{#2}\fi} \defaultdefine{Bbb}{\bf}
  \defaultdefine{frak}{\bf} \defaultdefine{=}{\B} 
  \defaultdefine{mathfrak}{\frak} \defaultdefine{mathbb}{\bf}
  \defaultdefine{mathcal}{\cal} \defaultdefine{implies}{\Rightarrow}
  \defaultdefine{beth}{BETH}\defaultdefine{cal}{\bf} \def\bbfI{{\Bbb I}}
  \def\mbox{\hbox} \def\text{\hbox} \def\om{\omega} \def\Cal#1{{\bf #1}}
  \def\pcf{pcf} \defaultdefine{cf}{cf} \defaultdefine{reals}{{\Bbb R}}
  \defaultdefine{real}{{\Bbb R}} \def\restriction{{|}} \def\club{CLUB}
  \def\w{\omega} \def\exist{\exists} \def\se{{\germ se}} \def\bb{{\bf b}}
  \def\equivalence{\equiv} \let\lt< \let\gt>
\providecommand{\bysame}{\leavevmode\hbox to3em{\hrulefill}\thinspace}
\providecommand{\MR}{\relax\ifhmode\unskip\space\fi MR }
\providecommand{\MRhref}[2]{%
  \href{http://www.ams.org/mathscinet-getitem?mr=#1}{#2}
}
\providecommand{\href}[2]{#2}
\begin{thebibliography}{}

\bibitem[Eis13]{Eis13}
Todd Eisworth, \emph{Getting more colors {II}}, J. Symbolic Logic \textbf{78}
  (2013), no.~1, 17--38.

\bibitem[Moo06]{Mo06}
Justin~Tatch Moore, \emph{{A solution to the $L$ space problem}}, Journal of
  the American Mathematical Society \textbf{19} (2006), 717--736.

\bibitem[Rin12]{Rin13}
Assaf Rinot, \emph{Transforming rectangles into squares, with applications to
  strong colorings}, Adv. Math. \textbf{231} (2012), no.~2, 1085--1099.

\bibitem[Tod87]{To2}
Stevo Todor\v{c}evi\'{c}, \emph{{Partitioning Pairs of Countable Ordinals}},
  Acta Math. \textbf{159} (1987), 261--294.

\bibitem[Sh:g]{Sh:g}
Saharon Shelah, \emph{{Cardinal Arithmetic}}, {Oxford Logic Guides}, vol.~29,
  {Oxford University Press}, 1994.

\bibitem[Sh:E12]{Sh:E12}
\bysame, \emph{{Analytical Guide and Corrections to \cite{Sh:g}.}},
  arxiv:math.LO/9906022.

\bibitem[Sh:276]{Sh:276}
\bysame, \emph{{Was Sierpi\'nski right? I}}, Israel Journal of Mathematics
  \textbf{62} (1988), 355--380.

\bibitem[Sh:280]{Sh:280}
\bysame, \emph{{Strong negative partition above the continuum}}, The Journal of
  Symbolic Logic \textbf{55} (1990), 21--31.

\bibitem[Sh:327]{Sh:327}
\bysame, \emph{{Strong negative partition relations below the continuum}}, Acta
  Mathematica Hungarica \textbf{58} (1991), 95--100.

\bibitem[Sh:572]{Sh:572}
\bysame, \emph{{Colouring and non-productivity of $\aleph_2$-cc}}, Annals of
  Pure and Applied Logic \textbf{84} (1997), 153--174, arxiv:math.LO/9609218.

\bibitem[JuSh:1025]{JuSh:1025}
Istvan Juhasz and Saharon Shelah, \emph{{Strong colorings yield
  $\kappa$-bounded spaces with discretely untouchable points}}, Proceedings of
  the AMS \textbf{143} (2015), 2241--2247, arxiv:math.LO/1307.1989.

\bibitem[Sh:F1422]{Sh:F1422}
Saharon Shelah, \emph{{Colouring of successor of regular, again}}.

\bibitem[Sh:F1783]{Sh:F1783}
\bysame, \emph{{Colourings for successor of regulars}}.

\end{thebibliography}

\end{document}